\documentclass[a4paper, 11pt]{article}

\usepackage{graphicx}
\usepackage{float}
\usepackage[utf8]{inputenc}
\usepackage{amsmath,amsthm,amssymb}
\usepackage{amsfonts}
\usepackage{color}
\usepackage{enumerate}
\usepackage{bm}
\usepackage{url}
\usepackage{subfigure}
\usepackage{color}

\newcommand{\R}{\mathbb{R}}
\newcommand{\Ss}{\mathbb{S}}
\newcommand{\M}{\mathbf{M}}
\newcommand{\NN}{\mathbf{N}}
\newcommand{\A}{\mathbf{A}}
\newcommand{\B}{\mathbf{B}}
\newcommand{\CC}{\mathbf{C}}
\newcommand{\X}{\mathbf{X}}
\newcommand{\Y}{\mathbf{Y}}
\newcommand{\x}{\mathbf{x}}
\newcommand{\w}{\mathbf{w}}
\newcommand{\cc}{\mathbf{c}}
\newcommand{\y}{\mathbf{y}}
\newcommand{\z}{\mathbf{z}}
\newcommand{\h}{\mathbf{h}}
\newcommand{\uu}{\mathbf{u}}
\newcommand{\vv}{\mathbf{v}}

\newcommand{\ee}{\mathbf{e}}

\newcommand{\D}{\mathbf{D}}
\newcommand{\T}{\mathbf{T}}
\newcommand{\J}{\mathbf{J}}
\newcommand{\Q}{\mathbf{Q}}
\newcommand{\I}{\mathbf{I}}
\newcommand{\G}{\mathrm{G}}
\newcommand{\N}{\mathrm{N}}
\newcommand{\Pp}{\mathcal{P}}
\newcommand{\Aa}{\mathcal{A}}
\newcommand{\EE}{\mathrm{E}}
\newcommand{\E}{\mathcal{E}}

\newcommand{\K}{\mathcal{K}}
\newcommand{\V}{\mathrm{V}}

\newcommand{\0}{\mathbf{0}}
\newcommand{\1}{\mathbf{1}}
\newcommand{\LL}{\mathbf{L}}
\newcommand{\Hh}{\mathbf{H}}
\newcommand{\Gg}{\mathbf{G}}

\newcommand{\rank}{\mathrm{rank}}

\newcommand{\diag}{\mathrm{diag}}
\newcommand{\col}{\mathrm{col}}
\newcommand{\vvec}{\mathrm{vec}}

\newtheorem{thm}{Theorem}

\newtheorem{lem}{Lemma}
\newtheorem{rem}{Remark}

\newtheorem{assum}{Assumption}

\begin{document}

\title{\bf Network Flows that Solve Sylvester Matrix Equations}
\author{Wen Deng\thanks{W. Deng is with LSC, Academy of Mathematics and Systems Science, Chinese Academy of Sciences; the School of Mathematical Sciences, University of Chinese Academy of Sciences, Beijing 100190, China;
		Australian Center for Field Robotics, School of Aerospace, Mechanical and Mechatronic Engineering, The University of Sydney, NSW 2006, Australia. (dengwena@amss.ac.cn).}~,~Yiguang Hong\thanks{Y.~Hong is with LSC, Academy of Mathematics and Systems Science, Chinese Academy of Sciences, Beijing 100190, China. (yghong@iss.ac.cn).}~,~Brian D. O. Anderson\thanks{B. D. O. Anderson is with Hangzhou Dianzi University, Hangzhou 310018, China; Data61-CSIRO, Australia; Research School of Electrical, Energy and Materials Engineering, The
		Australian National University, Canberra ACT 2600, Australia. (brian.anderson@anu.edu.au). B. D. O. Anderson is supported by the Australian Research Council (ARC) under grant DP-160104500.}~, and  Guodong Shi\thanks{G. Shi is with Australian Center for Field Robotics, School of Aerospace, Mechanical and Mechatronic Engineering, The University of Sydney, NSW 2006, Australia. (guodong.shi@sydney.edu.au).}
}
\date{}
\maketitle 
\begin{abstract}                  
	In this paper, we study distributed methods for solving a Sylvester equation in the form of  $ \A\X+\X\B=\CC $ for matrices $ \A, \B, \CC \in \R^{n\times n} $  with $ \X $ being the unknown variable.
	The entries of $ \A, \B $ and $ \CC $ (called data) are partitioned into a number of pieces (or sometimes we permit these pieces to overlap). 
	Then a network with a given structure is assigned, whose number of nodes is consistent with the partition. Each node has access to the corresponding set of data and holds a dynamic state. Nodes share their states among their neighbors defined from the network structure, and we aim to design flows that can asymptotically converge to a solution of this equation. 
	The decentralized data partitions may be resulted directly from networks consisting of physically isolated subsystems, or indirectly from artificial and strategic design for processing large data sets.
	 Natural partial row/column partitions, full row/column partitions and clustering block partitions of the data $ \A, \B $ and $ \CC $ are assisted by the use of the vectorized matrix equation. We show that the existing ``consensus $ + $ projection'' flow and the ``local conservation $ + $ global consensus'' flow for distributed linear algebraic equations can be used to drive distributed flows that solve this kind of equations. A ``consensus $ + $ projection $ + $ symmetrization'' flow is also developed for equations with symmetry constraints on the solution matrices. We  reveal some fundamental convergence rate limitations for such flows regardless of the choices of node interaction strengths and network structures. 
	 For a special case with $ \B=\A^T $, where the equation  mentioned is reduced to a classical Lyapunov equation, we demonstrate that by exploiting the symmetry of data, we can obtain  flows with lower complexity for certain partitions.
\end{abstract}

\textbf{Keywords:} Distributed algorithms; matrix equations; network flows; convergence rate.
\section{Introduction}
Recently, distributed optimization and computation in multi-agent networks have received growing research interest, where applications are witnessed in various problems for the control and operation of  large-scale network systems \cite{rabbat2004distributed,nedic2009distributed}.
A number of distributed  algorithms have arisen, involving many fields such as  distributed control and estimation, and distributed signal processing \cite{martinez2007motion,kar2012distributed,dimakis2010gossip}.
A related problem with growing research attention is to design
distributed algorithms for solving the linear algebraic equation $ \A\x=\mathbf{b} $ over a given network where the rows of $ \A $ and the entries of $ \mathbf{b} $ are allocated to individual nodes. These distributed optimization and computation ideas have also been explored in the areas of parallel computation and machine learning \cite{isard2007dryad,li2014scaling}, while efforts under the multi-agent frameworks focus more on scalability and resilience advantages for a given network structure.

As for the linear  equation $ \A\x=\mathbf{b} $, there are a few distributed solutions as discrete-time or continuous-time algorithms over networks \cite{lu2009distributed,lu2009distributed2,wang2014solving,anderson2015decentralized,shi2017networkTAC,wang2019scalable,liu2019arrow}. Every node only knows local information, such as one or several rows of $ \A $ and $ \mathbf{b} $, and then communicates with its neighbors about a dynamically evolving state. As long as  $ \A\x=\mathbf{b} $ has at least one solution, finding a solution to the original equation is equivalent to finding a solution in the intersection of affine subspaces defined by the solution spaces of individual nodes. With proper design of distributed flows, nodes can asymptotically agree on a certain solution to the overall equation $ \A\x=\mathbf{b}, $ complying with a given network structure and only exchanging state information (as opposed to information about $ \A $ and $ \mathbf{b} $). 
Notably, the ``consensus + projection'' flow \cite{shi2017networkTAC} has a simple form consisting of a standard consensus term and a local projection term onto every individual affine subspace. Generalized high-order flows with consensus and projection can even solve the equation approximately in the least-squares sense \cite{shi2017networkTAC}. In addition, a  double-layer network has been proposed to allow for a general data partition of the entries in $ \A $ and $ \mathbf{b} $, where the ``local conservation + global consensus'' flow \cite{wang2019scalable} or its variation can be used to solve $ \A\x=\mathbf{b} $ distributively.

Linear matrix equations, which are particular forms of  structured linear equations, appear in various fields of science and engineering \cite{zhou1996robust,bernstein2005matrix,simoncini2016computational}, such as the Sylvester equation in the form of $ \A\X+\X\B=\CC $  with $ \A\in\R^{n\times n}, \B\in\R^{m\times m}, \CC\in\R^{n\times m} $ and the unknown $ \X\in\R^{n\times m} $.  In fact, many Sylvester-type matrix equations in the control and automation areas serve as basic models for lots of fundamental systems and problems.
For example, the
Sylvester equation can be used to achieve pole/eigenstructure assignment by designing a controller for 
mechanical vibrating systems \cite{kim1999eigenstructure}, while
the Lyapunov equation with $ \B=\A^T $ plays an essential  role in studying the stability of linear time-invariant systems \cite{trentelman2012control}. 
The motivation for study distributed solver for Sylvester equations may come from the following two aspects: 
	 (\romannumeral1) Extension for linear algebraic equations to matrix equations is nontrivial, because the data partitions of  entries in $ (\A, \B, \CC) $ complying with a network would lead to fundamentally different computing problems compared to a standard linear algebraic equation;
	(\romannumeral2) Increasing growing study of complex network systems requires distributed solutions of the matrix equations from physically isolated data sets to problems as basic as stability validation.

In this paper, we concentrate on seeking a solution to the matrix equation $ \A\X+\X\B=\CC $ with $ \A, \B, \CC, \X \in \R^{n\times n} $ in a distributed way.
Here we choose to work on square matrices to facilitate a simplified presentation; nonetheless, our methods and analysis can be straightforwardly generalized to the general Sylvester equation with $ \X\in\R^{n\times m} $, because $ \X $ being a square matrix plays no role in our algorithm design and convergence characterizations. Note that the system cannot be directly studied with the methods for the equation $ \A\X\B=\mathbf{F} $ discussed in \cite{zeng2018distributedlmeTAC} because these two equations necessarily give rise to  different patterns of assignment of data to network nodes. The work \cite{zeng2018distributedlmeTAC} builds a solution procedure from an optimization perspective and solves several primal and dual optimization problems via distributed methods, while we plan to transform  matrix equations by vectorization and take  advantage of the above referenced distributed algorithms for solving  $ \A\x=\mathbf{b} $. 
More concretely, in our design, each node has access to local data in matrices $ \A, \B $ and $ \CC $ with the following several partition patterns,  which may be suited to certain different problems.
\begin{itemize}
	\item [(\romannumeral1)]
	[Partial  row/column partition] E.g., for an $ n $-node network, each node $ i $ holds the $ i $-th column of $ \B $ and $ \CC $, with the entire $ \A $ known to the whole network. We show that with such partition,
	we can utilize the ``consensus $ + $ projection''  flow for an $ n $-node network, under which we establish convergence with an explicit rate and more interestingly, $  $a rate limitation characterization of the flow. In addition, we design a ``consensus $ + $ projection $ + $ symmetrization'' flow for symmetry constraints on the solution matrices, followed by its  properties of convergence and rate limitation.
    \item [(\romannumeral2)]
    [Full  row/column partition] E.g., for an $ n $-node network, each node $ i $ holds the $ i $-th row of $ \A $, and the $ i $-th column of $ \B $ and $ \CC $.  We show that under this type of partition, the Sylvester equation can be solved distributedly by introducing an auxiliary variable and taking advantage of an augmented ``consensus + projection'' flow in a node state space with dimension $ n^2(1+n) $.
	\item [(\romannumeral3)]
	[Clustering block partition] E.g., for a double-layer network with $ n $ clusters, the $ i $-th of which having $ n $ nodes holds the entire $ \A $, and the $ i $-th column of $ \B $ and 
	$ \CC $, where each node $ j $ within cluster $ i $ is assigned to the $ (j,i) $-th entry of $ \B $ and $ \CC $, and additional matrix $ \A $ if $ j=i $. Taking advantage of the ``local conservation $ + $ global consensus'' flow, we establish convergence with an explicit rate as well. As a byproduct of the study, a fundamental property of the convergence rate  in the ``local conservation $ + $ global consensus'' flow is also established, which is of independent interest.
\end{itemize}
In a brief discussion, we also show that the data $ \A, \B $ and $ \CC $ can be partitioned over an $ n^2 $-node network, where each node  holds one row of $ \A $, one column of $ \B $ and one entry of $ \CC $. As a result, the data complexity at each node is reduced with $ n^2 $ nodes, while the rate of convergence for the resulting flow, however, becomes lower due to the increased network size.
If in addition, there is a particular case with $ \B=\A^T $, where the equation becomes a Lyapunov equation, and  by exploiting the symmetry, the size of nodes can be reduced to $ n(n+1)/2 $ compared with the $ n^2 $-node network.

For this paper, the remainder is organized as follows. In Section \ref{sec:problem def}, we define the considered matrix equation problem with a motivating example. In Section \ref{sec:row/column partition}, we present a network flow with partial row/column partitions and prove the convergence rate limitation, followed by some numerical examples and discussions. We consider full  row/column partitions in Section \ref{sec:row+column partition} with corresponding network flow. In Section \ref{sec:double partition}, we present a network flows with a clustering block partition and set out some properties. In Section \ref{sec:conclusion}, we conclude the paper briefly with a few remarks. Finally, some useful results related to linear algebra, projection, and exponential stability are introduced in Appendix \ref{sec:preliminaries}, and other details of proofs are given in subsequent appendices.

{\em Notation}: Let $ \0 $ or $ \1 $  represent the matrix (or vector) with all entries being $ 0 $ or $ 1 $, and their dimensions are indicated by subscripts. Let $ \I_n $ denote an $ n $ by $ n $ identity matrix and $ \ee_i $ denote the $ i $-th column vector of $ \I_n. $ 
  Let $ \col\{\M_{[1]},\cdots, \M_{[n]}\}=[\M_{[1]}^T,\cdots, \M_{[n]}^T]^T $ be a stack of matrices $ \M_{[i]}, i=1,\cdots,n $. Let 
   $ \vvec_{mn}(\cdot) $ be a mapping from $ \R^{m\times n} $ to $ \R^{mn} $: $ \vvec(\A)=\col\{\A_1,\cdots,\A_n\} $ with $ \A_i $ being the $ i $-th column of $ \A. $ The inverse mapping of $ \vvec_{mn}(\cdot) $ can be well-defined, which is denoted by  $ \vvec_{mn}^{-1}(\cdot) $. The subscripts of $ \vvec_{mn}(\cdot) $ and $ \vvec_{mn}^{-1}(\cdot) $ would be dropped  whenever there is no ambiguity of the space dimensions.
  Denote by $ \mathrm{span}(\M), \ker(\M) $ and $ \rank(\M) $, the column space, the kernel space and the rank of a matrix $ \M $, respectively. Let $\mathrm{diag}\{\mathbf{F}_{[1]},\cdots,\mathbf{F}_{[n]}\} $ denote the block diagonal matrix with sub-blocks $ \mathbf{F}_{[i]}, i=1,\cdots,n $. For a matrix $ \A\in\R^{n\times n} $ with all real eigenvalues, let $ \mathrm{spec}(\A)=\{\lambda_1(\A),\cdots,\lambda_n(\A)\} $ denote the set of all the eigenvalues of $ \A $ with $ \lambda_1(\A)\ge \lambda_2(\A)\ge \cdots\ge \lambda_{n}(\A) $.
Let $ \otimes $ denote the Kronecker product and $ \mathrm{dim}(\mathcal{V}) $ represent the dimension of a subspace $ \mathcal{V} $ in $ \R^n. $ Let $ \|\cdot\| ( \|\cdot\|_F)$ denote the Euclidean (Frobenius) norm of a vector (matrix) and $ B_{\delta}:=\{\x\in\R^n: \|\x\|\le\delta\} $ denote the closed ball with radius $ \delta $ and  center at the origin.
Denote a network graph $ \G=(\V,\EE) $ with node set $ \V $ and edge set $ \EE. $ All graphs  in the remainder of this paper  are  connected and undirected.
The neighbor set of node $ i $ is given by $ \N_i:=\{j: (i,j)\in \EE\} $, from which the node $ i $ can receive  information. Introduce $ \A_{\G} $ as the adjacency matrix of $ \G, $ with $ [\A_{\G}]_{ij}=1 $ if $ (i,j)\in \EE, $ and $ [\A_{\G}]_{ij}=0, $ otherwise.
The Laplacian matrix $ \LL_{\G} $ associated with $ \G $  is defined by $ \LL_{\G}=\D_{\G}-\A_{\G}, $ where $ \D_{\G}=\mathrm{diag}\{\sum_{j\in \N_i}[\A_{\G}]_{ij}, i\in\V\} $ is the degree matrix of the graph $ \G $.

\section{Problem Definition} \label{sec:problem def}
In this section, we introduce the motivation of the study for matrix equations over networks and define the problem of interest.
\subsection{Matrix Equation}
Consider a matrix equation with respect to variable $ \X\in \R^{n\times n} $:
\begin{equation}\label{eq:Sylvester}
\A\X+\X\B=\CC, \quad \A, \B, \CC\in\R^{n\times n}.
\end{equation}
By vectorization, we have the following equivalent equation of \eqref{eq:Sylvester}
\begin{equation}\label{eq:olinear_equation}
(\I_n\otimes \A+\B^T\otimes \I_n)\x=\cc,
\end{equation}
where $ \x=\mathrm{vec}(\X), \cc=\mathrm{vec}(\CC) $. 
There are three cases covering the solvability properties.
\begin{enumerate}
	\item [(\uppercase\expandafter{\romannumeral1})] The solution to $ \eqref{eq:olinear_equation} $ is unique if and only if the matrix $ \I_n\otimes \A+\B^T\otimes \I_n $ is nonsingular, which is equivalent to $ \mathrm{spec}(\A)\cap\mathrm{spec}(-\B)=\emptyset $  \cite{horn1990matrix}.
	\item [(\uppercase\expandafter{\romannumeral2})] The solution to $ \eqref{eq:olinear_equation} $ is an infinite set when $ \mathrm{spec}(\A)\cap\mathrm{spec}(-\B)\ne\emptyset $ and $ \cc\in \mathrm{span}(\I_n\otimes \A+\B^T\otimes \I_n). $
	\item [(\uppercase\expandafter{\romannumeral3})] There is no exact solutions to \eqref{eq:olinear_equation}
	when $ \cc\notin \mathrm{span}(\I_n\otimes \A+\B^T\otimes \I_n). $
\end{enumerate}
\subsection{A Motivating Example}\label{subsec:motivating-ex}

Consider the following network system with $ n $ dynamically coupled subsystems, for $ i\in\V=\{1,\cdots,n\} $:
\begin{equation}\label{eq:subsys}
\dot{\y}_i=\D_{ii}\y_i+\sum_{j=1,j\ne i}^n\D_{ij}\y_j, \quad \D_{ij}\in\R^{m\times m}, \y_i\in\R^m, 
\end{equation} 
where $ \y_i $ is the state of the subsystem $ i $ and $ \D_{ij} $ represents the dynamical influence from subsystem $ j $ to subsystem  $ i $. The system \eqref{eq:subsys} is arguably one of the most basic models for dynamical networks with linear couplings, which may represent a large number of practical network systems ranging from power distribution, transportation, and controlled formation \cite{fuhrmann2015mathematics,trumpf2019controllability,fax2004information,oh2015survey,blackhall2010structural}.
The overall network dynamics is in the form of $ \dot{\y}(t)=\A\y(t) $, where $ \y(t)=\col\{\y_{1}(t),\cdots,\y_n(t)\} $ is the network state and 
\begin{equation}
\A=\begin{bmatrix}
\D_{11} & \D_{12} &\cdots & \D_{1n}\\
\D_{21}& \D_{22} &  \cdots & \D_{2n}\\
\vdots & \vdots & \cdots & \vdots\\
\D_{n1}& \D_{n2}&\cdots & \D_{nn}
\end{bmatrix}.
\end{equation}

We introduce the following problem.

\noindent\textbf{Problem:} Each subsystem $ i $ knows  $ \D_{ij}, j\in\V, $ and aims to verify the stability of the overall network system in a distributed manner without directly revealing its dynamics $ \D_{ij} $ to any other nodes. 

Here, the words ``in a distributed manner'' imply that the subsystem $ i $ only interacts with a set of neighbors over a communication graph $ \G=(\V,\EE) $. The communication graph $ \G $ may or may not coincide with interaction graph encoded in the dynamics \eqref{eq:subsys}: $ \G_{\mathbf{L}}(\V,\EE_{\mathbf{L}}) $ with $ (j,i)\in\EE_{\mathbf{L}} $ if and only if $ \D_{ij}\ne0. $
If   the $ i $-th subsystem can hold a dynamical state $ \X_i\in\R^{d\times d} $, which is shared over the communicating links over the graph $ \G, $  then any of the subsystems can verify the stability of the overall network if $ \X_i(t) $ converges to a positive definite solution to the following Lyapunov equation (\cite{khalil2002nonlinear}):
\begin{equation}\label{eq:first-moexam-lyap}
\A\mathbf{X}+\mathbf{X}\A^T=-\I_{d} ,\quad \A,\X\in\R^{d\times d},\ d=mn.
\end{equation}
Therefore, distributed solvers of the Sylvester matrix equations may be used as a tool for stability validation of network systems.

\subsection{Problem of Interest}

We impose the following assumption, which holds throughout the rest of the article.
\begin{assum}\label{assumption-exact solution}
	Equation \eqref{eq:Sylvester} has at least one exact solution.
\end{assum}

Under Assumption \ref{assumption-exact solution}, we focus on solving the matrix equation $ \A\X+\X\B=\CC $ with solution case (\uppercase\expandafter{\romannumeral1}) or (\uppercase\expandafter{\romannumeral2}) in a distributed manner. To be precise, we mainly aim to 
\begin{itemize}
	\item [(\romannumeral1)] distribute the entries of $ \A, \B $ and $ \CC $ over the nodes in a network $ \G=(\V,\EE) $; 
	\item [(\romannumeral2)] assign each node a dynamic state which can be shared among the neighbors over $ \G $; 
	\item [(\romannumeral3)] design decentralized flows that drive the states of nodes to the solutions of the Sylvester equation;
	\item [(\romannumeral4)] explore the convergence and the limitation of the convergence rate.
\end{itemize}

In our motivating example on stability validation of network systems, the data partition is due to the natural isolation of subsystems.
 The advantage of data partition also arises from the fact that a large data set with the size of $ 3n^2 $ can be partially split into  multiple subsets of reduced size and  handled in a distributed way. Similar ideas have been explored in distributed convex optimization \cite{nedic2010constrainedTAC} and submodular optimization \cite{mirzasoleiman2016distributed}.

\section{Partial Row/Column Partition}\label{sec:row/column partition}
In this section, we consider the data partitions where the entire $ \A $ with  partial $ \B $ and $ \CC $, or the entire $ \B $ with  partial $ \A $ and $ \CC $  would be allocated at $ n $ individual nodes.

\subsection{Partial Column/Row Partition}\label{subsec:row/column part}
Denote $ [\A_0]_i:=[\0,\cdots,\A,\cdots,\0], $ where the  $ i $-th block is $ \A $ and the other $ n-1 $ blocks are $ \0_{n\times n} $, and $ [\B_0^T]_i:=[\0,\cdots,\B^T,\cdots,\0] $ as well.
Over an $ n $-node network, we consider two main partitions as follows.
\begin{enumerate}
	\item[(\romannumeral1)] $ [ \B$-$\CC $ Column Partition]   Node $ i $ holds $ \A, $ and the $ i $-th column of $ \B $ and  $ \CC $, denoted by $ \B_i $ and $ \CC_i $, respectively. Equivalently, node $ i $ has access to an equation
	\begin{equation}\label{eq:decomposition-1}
	([\A_0]_i+\B_i^T\otimes \I_n)\mathrm{vec}(\X)=\CC_i.
	\end{equation}
	\item[(\romannumeral2)] $ [ \A$-$\CC $ Row Partition] Node $ i $ holds $ \B, $ and the $ i $-th row of $ \A $ and  $ \CC $, denoted by $(\A^T)_i^T $ and $ (\CC^T)_i^T $, respectively. Equivalently, node $ i $ has access to an equation
	\begin{equation}\label{eq:decomposition-1'}
	([\B_0^T]_i+(\A^T)_i^T\otimes \I_n)\mathrm{vec}(\X^T)=(\CC^T)_i.
	\end{equation}
\end{enumerate}
\begin{rem}
	Except for the two partitions above, there may be  other partitions, such as
	 the entire $ \A $ with $ \B $ Column/$ \CC $ Row (or $ \B $ Row/$ \CC $ Column, or $ \B$-$\CC $ Row) Partition.
	It turns out those partitions will have a  different nature and be suitable for different algorithms. Nevertheless, due to the feature of  the partitions $\B$-$\CC $ Column and $\A$-$\CC $ Row, the matrix equation \eqref{eq:Sylvester} can be easily reformulated into \eqref{eq:decomposition-1} and \eqref{eq:decomposition-1'}, which are concisely shown as $ n $ separate linear algebraic equations. Therefore, the  $ \B$-$\CC $ Column Partition and  $ \A$-$\CC $ Row Partition are suitable for the ``consensus $ + $ projection'' flow.
\end{rem}

  In fact, these two
 partitions  are essentially equivalent from an algorithmic point of view because $ \A\X+\X\B=\CC $ is equivalent to $ \B^T\X^T+\X^T\A^T=\CC^T. $ Therefore, in the following we  focus on the $\B$-$\CC $ Column Partition.
Define
\[ \E_i:=\{\y\in \R^{n^2}:([\A_0]_i+\B_i^T\otimes \I_n)\y=\CC_i\},\quad i\in\V,\]
 where $ \E_i $ is an affine subspace and  $ \V=\{1,\cdots,n\}  $.
It follows from the solution cases where Case  (\uppercase\expandafter{\romannumeral1})  means $ \E:=\cap_{i=1}^n\E_i $ is a singleton; Case  (\uppercase\expandafter{\romannumeral2})  means $ \E $ is an affine space with a nontrivial dimension; and
Case  (\uppercase\expandafter{\romannumeral3})  means $ \E=\emptyset $.

\begin{rem}
	We have assumed that there are $ n $ nodes with node $ i $ having partial data $  \A, \B_i $ and $ \CC_i. $ We could if desired assume that  the number of nodes $ p $ is less than $ n $; then we  have a partition where node $ i $ has access to $ \A, \{\B_{i(1)},\cdots,\B_{i(q_i)}\} $ and $ \{\CC_{i(1)},\cdots,\CC_{i(q_i)}\} $ for some $ q_i<n $ with $ i=1,\cdots,p $ and $ i(r)\in \{1,\cdots,n\} $. In this scenario, we readjust the affine subspace $ \E_i $ to
	$ \E_i:=\{\y\in \R^{n^2}:([\A_0]_k+\B_k^T\otimes \I_n)\y=\CC_k, k=i(1),\cdots,i(q_i)\}, $
	where the index of $ i $ satisfies $
	\cup_{i=1}^p\{i(1),\cdots,i(q_i)\}=\V. $ Case  (\uppercase\expandafter{\romannumeral1}) and Case  (\uppercase\expandafter{\romannumeral2}) can guarantee that every $ \E_i $ and the 
	intersection are nonempty. Our discussion can be applied to this generalized partition, 
	and the determination of the best way to form the partition, from the viewpoint of convergence rate or communications burden, etc., should be under consideration according to  specific circumstances.
\end{rem}

\subsection{Generalized ``Consensus $ + $ Projection'' Flow}
Let a mapping $ \Pp_{\E_i}: \R^{n^2}\to \R^{n^2} $ be the projector onto the affine subspace $ \E_i $ and $ K>0 $ be a given constant.
Motivated by references \cite{shi2017networkTAC,shi2013reachingTAC}, we consider the following continuous-time network flow:
\begin{equation}\label{eq:Flow-n-nodes}
\dot{\x}_i=K\Big(\sum_{j\in \N_i}(\x_j-\x_i)\Big)+\Pp_{\E_i}(\x_i)-\x_i, \quad i\in\V,
\end{equation}
where $ \x_{i}\in \R^{n^2} $ is a state held by node $ i $. Note that we could, if desired, insert a further multiplication $ K^{'} $ say of the term $ \Pp_{\E_i}(\x_i)-\x_i. $ This can be expected to change the convergence rate up and down. Obviously also, if $ K $ and $ K^{'} $ were both to be scaled by the same amount, the convergence rate can be changed. To separate these two effects, in this paper we select $ K^{'}=1 $, and consider the effect of adjusting $ K $ alone. 
The flow \eqref{eq:Flow-n-nodes} is the so-called ``Consensus + Projection'' Flow proposed in \cite{shi2017networkTAC}, where for the problem under consideration each $ \E_i $ is  an  affine subspace of $ \R^{n^2} $ 
as considered originally in \cite{shi2017networkTAC}.

Define $ 
\Hh_i:=
[\A_0]_i+\B_i^T\otimes \I_n,  i\in\V, $ and
 $ \Hh:=\col\{\Hh_{1},\cdots,\Hh_{n}\}=\I_n\otimes \A+\B^T\otimes \I_n. $ 
In view of Lemma \ref{lemma:projection}, the flow \eqref{eq:Flow-n-nodes} can be written in a compact form for $ \x=\col\{{\x_1},\cdots,{\x_n}\}\in \R^{n^3}, $
\begin{equation}\label{eq:Flow-X}
\dot{\x}=-(K\LL_{\G}\otimes \I_{n^2}+\J)\x+\Q_{\CC},
\end{equation}
where $ \LL_{\G} $ is the Laplacian matrix, $ \J $ is a block-diagonal matrix 
$ \mathrm{diag}\{\Hh^{\dag}_1\Hh_1,\cdots,\Hh^{\dag}_n\Hh_n\}\in \mathbb{R}^{n^3\times n^3} $ with $ \Hh^{\dag}_i $ being a M-P pseudoinverse of $ \Hh_i $, and $ \Q_{\CC}:=\col\{\Hh^{\dag}_1\CC_1, \cdots, \Hh^{\dag}_n\CC_n\}. $
The existence of an equilibrium point of system \eqref{eq:Flow-X} is guaranteed by Assumption \ref{assumption-exact solution}. In fact, if $ \uu_0\in \cap_{i=1}^n\{\y:\Hh_i\y=\CC_i\}, $ let $ \uu^*=\col\{\uu_1^{*},\cdots,\uu_n^{*}\}=\1_n\otimes\uu_0\in \R^{n^3}. $ We have $ \Hh_i\uu_i^*-\CC_i=\0, $ furthermore, $ \Hh^{\dag}_i\Hh_i\uu_i^*-\Hh^{\dag}_i\CC_i=\0 $ (namely, 
$ \J\uu^*-\Q_{\CC}=\0 $).
Combining with $ (K\LL_{\G}\otimes \I_{n^2})\uu^*=\0, $ we conclude that
$ \uu^* $ is an equilibrium of \eqref{eq:Flow-X}. In the event that equation \eqref{eq:Sylvester} has a unique solution, $ \uu^* $ is also unique, an almost immediate consequence of the following theorem.

Denote $ \J_{\LL}:=K\LL_{\G}\otimes \I_{n^2}+\J $ and  $ r(K)=\min\{\lambda\in \mathrm{spec}(\J_{\LL}),\lambda\ne0\}. $
Recall that $ \lambda_{k}(\M) $ represents the $ k $-th largest eigenvalue of a symmetric matrix $ \M $. 
The flow \eqref{eq:Flow-n-nodes} has a fundamental convergence rate limitation established precisely in Theorem \ref{thm:convergence}.

\begin{thm}\label{thm:convergence}
	Under the $\B$-$\CC $ Column Partition, for any initial value $ \x
	_0=\col\{\x_1(0),\cdots,\x_n(0)\} $, there exists $ \X^*(\x_0)\in\R^{n\times n} $ as a solution to \eqref{eq:Sylvester}, such that along the flow \eqref{eq:Flow-n-nodes} $ \vvec^{-1}(\x_{i}(t)) $ converges to $ \X^*(\x_0) $ exponentially, for all $ i\in \V$. To be precise, the following statements hold.
	\begin{enumerate} 
		\item[(\romannumeral1)] For any $ i\in \V, $ \begin{align*} \lim\limits_{t\to \infty}\vvec^{-1}(\x_{i}(t))=&\X^*(\x_0)=\vvec^{-1}\Big(\frac{1}{n}\sum_{i=1}^{n}\Pp_{\cap_{i=1}^n\E_i}(\x_i(0))\Big). \end{align*}
		\item[(\romannumeral2)] There exist $ \beta(\x_0), r(K)>0 $, such that for all $ t\ge 0 $,
		$$  \sum_{i=1}^n\|\vvec^{-1}(\x_{i}(t))- \X^*(\x_0)\|_F^2\le \beta(\x_0) e^{-2r(K)t}, $$
		where the exponential rate $ r(K) $ is a non-decreasing and bounded  function with respect to $ K $ satisfying $ \lim\limits_{K\to \infty}r(K)= \lambda_{\rank(\Hh)}(\frac{1}{n}(\sum_{i=1}^{n}\Hh^{\dag}_i\Hh_i)) $.
	\end{enumerate}
\end{thm}
Details of the proof for Theorem \ref{thm:convergence} can be found in Appendix \ref{App:1con}.
\begin{rem}
	Define \begin{equation*}
		\begin{array}{l}
	\lambda^{\ast}=\max\{\lambda_1(\Hh_{1}\Hh_{1}^T),\lambda_1(\Hh_{2}\Hh_{2}^T),\cdots,\lambda_1(\Hh_{n}\Hh_{n}^T)\}, \\ \lambda_{\ast}=\min\{\lambda_n(\Hh_{1}\Hh_{1}^T),\lambda_n(\Hh_{2}\Hh_{2}^T),\cdots,\lambda_n(\Hh_{n}\Hh_{n}^T)\},\\
	f_{\A\B}=\I_n\otimes(\A^T\A)+(\sum_{i=1}^n\B_i\B_i^T)\otimes\I_{n}+\B\otimes\A+(\B\otimes\A)^T.
	\end{array}\end{equation*} Then,  
	if every $ \Hh_i=
	[\A_0]_i+\B_i^T\otimes \I_n $ has  full row rank, i.e. $ \rank(\Hh_i)=n, $ we have
	\begin{align*}
	\lambda_{n^2}\left(\frac{f_{\A\B}}{n\lambda^{\ast}}\right)\le\lim\limits_{K\to \infty}r(K)\le\lambda_{1}\left(\frac{f_{\A\B}}{n\lambda_{\ast}}\right).
	\end{align*} 
\end{rem}

\begin{rem}\label{rem:lscase}
	Given data $ \A, \B $ and $ \CC $, it might be tedious if not impossible to verify the solvability conditions of Case (\uppercase\expandafter{\romannumeral1}) and (\uppercase\expandafter{\romannumeral2}) (as defined according to \eqref{eq:olinear_equation}), or it might be that the data corresponds to the Case (\uppercase\expandafter{\romannumeral3}). Hence, we could consider the least-squares solution in the sense of $ \min \sum_{i=1}^n\|\Hh_i\x_i-\CC_i\|^2  $ using similar ideas. Inspired by \cite{shi2017networkTAC}, when $ \Hh $ has full column rank, we can use the flow 
	\begin{equation}\label{eq:Flow-n-nodes-ls}
		\dot{\x}_i=K\big(\sum_{j\in \N_i}(\x_j-\x_i)\big)-\Hh_i^{\dag}(\Hh_i\x_i-\CC_{i}), \quad i\in\V.
		\end{equation}
	 For any $ \epsilon>0 $, there exists  $ K_0(\epsilon)>0 $, such that every $ \x_i(t) $ converges to the $ \epsilon $-neighborhood of the least-squares solution (e.g.,  Theorem 6 in \cite{shi2017networkTAC}) if $ K\ge K_0(\epsilon) $. 
\end{rem}

\subsection{``Consensus $ + $ Projection $ + $ Symmetrization'' Flow}
It would also be of interest to find a symmetric solution to \eqref{eq:Sylvester} if indeed \eqref{eq:Sylvester} admits at least one symmetric solution. The ``consensus $ + $ projection'' flow \eqref{eq:Flow-n-nodes} however cannot guarantee to find such a symmetric solution.
Let a mapping $ \Pp_{\Ss_{nn}}:\R^{n^2}\to\R^{n^2} $ be the projector onto a subspace $$ \Ss_{nn}:=\{\y\in\R^{n^2}:\y=\vvec(\X), \ \mathrm{for\ some\ symmetric\ matrix\ } \X\in\R^{n\times n}\}  $$  and $ K, K_s>0 $ be given constants.
We propose the following ``consensus $ + $ projection $ + $ symmetrization'' flow, for $ i\in \V $:
\begin{equation}\label{eq:Flow-n-nodes and symmetric}
\begin{aligned}
\dot{\x}_i=K\Big(\sum_{j\in \N_i}(\x_j-\x_i)\Big)+\Pp_{\E_i}(\x_i)-\x_i+K_s\Big(\Pp_{\Ss_{nn}}(\x_i)-\x_i\Big).
\end{aligned}
\end{equation} 
The additional term $ K_s(\Pp_{\Ss_{nn}}(\x_i)-\x_i) $  plays a role in driving the node states to $ \Ss_{nn} $.
 Then we present the following result.
\begin{thm}\label{thm:sym-convergence}
	Suppose that there is  a symmetric  solution to \eqref{eq:Sylvester}. Then, under the $\B$-$\CC $ Column Partition, for any initial value $ \x_0=\col\{\x_1(0),\cdots,\x_n(0)\} $, there exists $ \X_s^*(\x_0)\in\R^{n\times n} $ as a symmetric solution to \eqref{eq:Sylvester}, such that  along the flow \eqref{eq:Flow-n-nodes and symmetric}  $ \vvec^{-1}(\x_i(t)) $ converges to $ \X_s^*(\x_0) $ exponentially, for all $ i\in\V $. Moreover, there
exist $ \beta_s(\x_0), r_s(K,K_s)>0 $, such that for all  $ t\ge0, $
		$$  \sum_{i=1}^n\|\vvec^{-1}(\x_{i}(t))- \X_s^*(\x_0)\|_F\le \beta_s(\x_0) e^{-2r_s(K,K_s)t}. $$ In fact, the   exponential rate $ r_s(K,K_s) $ satisfies $$ r_s(K,K_s)\le\min\{1+K_s, 1+K\lambda_1(\LL_{\G})\}, $$ for all $ K, K_s>0 $, where $ \LL_{\G} $ is the Laplacian matrix of the relevant graph.		 	
\end{thm}
The proof of Theorem \ref{thm:sym-convergence} is in Appendix  \ref{App:sym-con}.

\subsection{Numerical Examples} \label{sec:example-n-nodes}
In this part, we present several numerical examples.

\noindent{\bf Example 1. }
\begin{figure}[htbp]	
		\centering
		\includegraphics[width=0.4\linewidth]{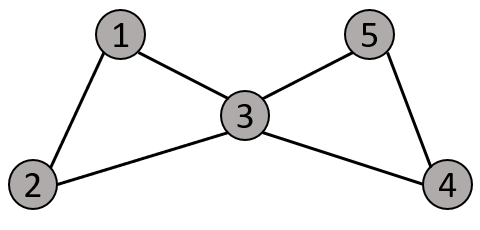}
		\caption{The graph structure of  five nodes.}
		\label{fig:5-nodes}	
\end{figure}
Consider a matrix equation:
\begin{equation}\label{ex1:eq}
\A\X+\X\B=\CC, \A, \B, \CC\in\R^{5\times 5},
\end{equation}  
where \begin{align*}
&\A=[7\	1\	1\	1\	5;7\	2\	8\	3\	0;1\	4\	8\	7\	7;7\	8\	4\	6\	7;5\	8\	6\	8\	5];\\
&\B=[6\	6\	7\	4\	4;6\	0\	6\	3\	4;3\	2\	3\	6\	5;5\	0\	8\	6\	6;1\	1\	0\	1\	6];\\
&\CC=[2\	 4\	6\	8\	7;5\	8\	2\	4\	2;5\	3\	4\	1\	7;1\	5\	6\	1\	2;1\	2\	7\	2\	7].
\end{align*}
It can be verified that this equation has a unique solution $ \X^* $. 
The related $ 5 $ nodes in an interconnected network  forms a  graph  shown in Fig. \ref{fig:5-nodes}.
\begin{figure}[h!]
	\centering
	\includegraphics[width=0.5\linewidth]{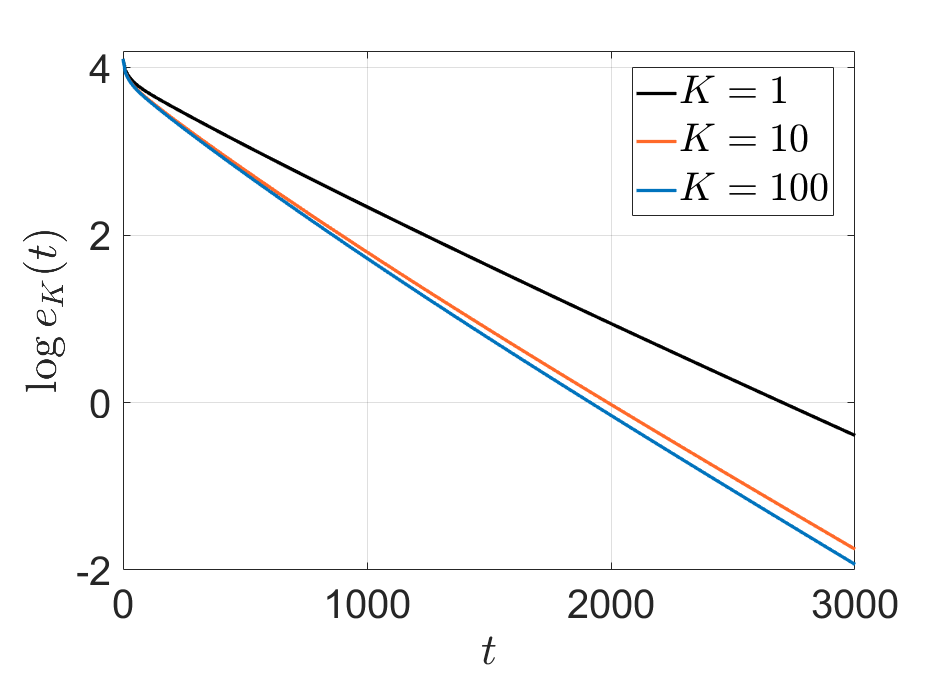}
	\caption{The evolution of $ \log e_K(t) $ for $ K=1,10,100, $ respectively.}
	\label{fig:k110100}	
\end{figure}
\begin{figure}[htbp]
	\centering
	\includegraphics[width=0.5\linewidth]{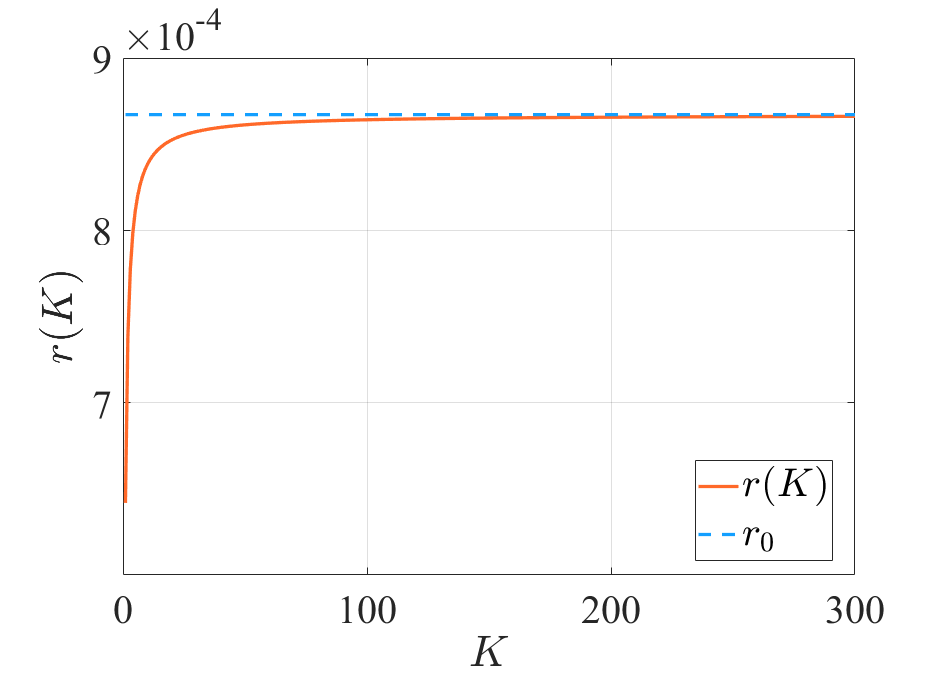}
	\caption{The trajectories of $ r(K) $ over $ K $ and the reference $ r_0 $.}
	\label{fig:rKt}
\end{figure}
Taking the initial value to be $ \0 $, we plot the trajectories of 
\begin{equation}\label{eq:et}
e_K(t):=\sum_{i=1}^5\|\x_i(t)-\mathrm{vec}(\X^*)\|^2,
\end{equation}
in logarithmic scales 
for $ \x_i(t) $ evolving along \eqref{eq:Flow-n-nodes} with $ K=1,10,100, $ respectively, in Fig. \ref{fig:k110100}, which  validates the exponential convergence in Theorem \ref{thm:convergence}. With different values of  $ K, $ we  calculate $ r(K) $ and plot $ r(K) $ over $ K $ in Fig. \ref{fig:rKt} with $ r_0=\lambda_{\rank(\Hh)}((\sum_{i=1}^{n}\Hh^{\dag}_i\Hh_i)/n) $  drawn as a reference. Fig. \ref{fig:rKt} shows that $ r(K) $ increases as  $ K $ increases and  $ r(K) $ always has an upper bound, which is consistent with $ r_0=\lim\limits_{K\to\infty}r(K) $.

\noindent{\bf Example 2. }
Still consider an equation in the form of \eqref{ex1:eq} and an interconnected network shown in Fig. \ref{fig:5-nodes}. We investigate  two sets of data for $ \A, \B, \CC $ as
\begin{align*}
&\A^{[1]}=[0\	0\	0\	5\	0;
0\	2\	0\	0\	2;
1\	3\	0\	0\	0;
0\	0\	4\	0\	0;
0\	0\	0\	0\	0];\\
&\B^{[1]}=[7\	4\	4\	7\	10;
3\	8\	6\	7\	3;
10\	8\	7\	2\	6;
0\	2\	8\	1\	2;
4\	5\	3\	5\	8];\\
&\CC^{[1]}=[8\	1\	6\	8\	3;
8\	5\	5\	3\	7;
4\	8\	0\	5\	7;
6\	9\	3\	2\	7;
1\	1\	2\	6\	5]; \\
&\A^{[2]}=[1\	5\	10\	1\	9;
2\	10\	0\	4\	2;
9\	1\	8\	3\	3;
2\	4\	8\	8\	1;
8\	1\	9\	4\	1];\\
&\B^{[2]}=[0\	0\	8\	6\	0;
2\	0\	0\	0\	0;
0\	1\	0\	0\	0;
0\	0\	0\	0\	0;
0\	0\	0\	0\	3];\\
&\CC^{[2]}=[9\	6\	2\	0\	3;
6\	4\	1\	9\	9;
5\	5\	2\	9\	4;
1\	4\	2\	5\	1;
9\	1\	4\	5\	8].
\end{align*}	
It can be seen that (\romannumeral1) $ \A^{[1]} $ is sparse, $ \B^{[1]} $ is dense; (\romannumeral2) $ \A^{[2]} $ is dense, $ \B^{[2]} $ is sparse.  
\begin{figure}[htbp]
	\begin{minipage}{0.495\linewidth}
		\centering
		\includegraphics[width=1\linewidth]{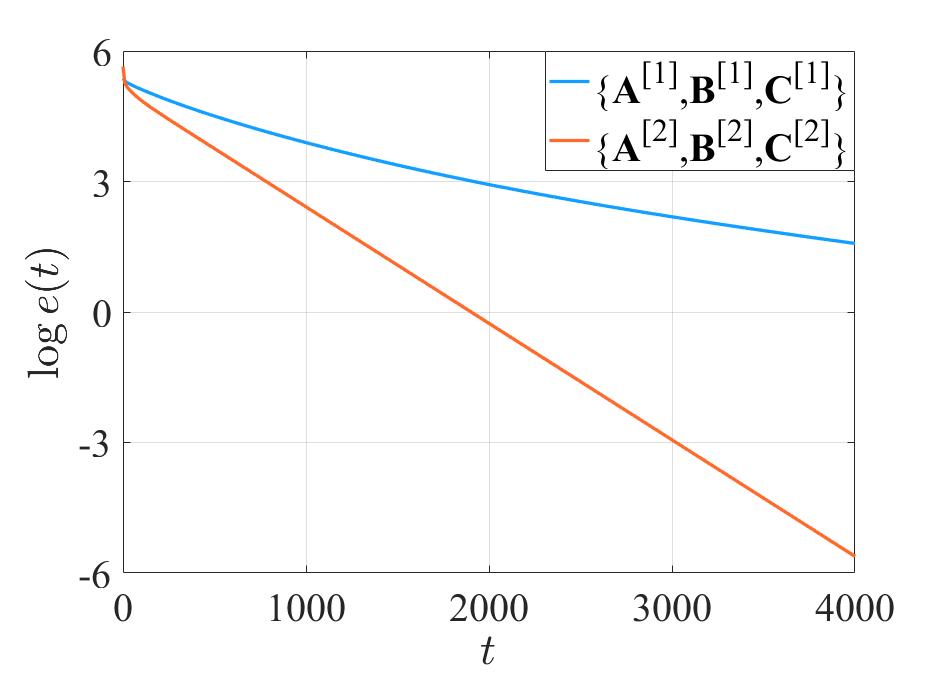}
		\caption{The trajectories of $ \log e(t) $ for two problems with data sets  $ \{\A^{[i]},\B^{[i]},\CC^{[i]}\}, i=1,2 $ under $ \B$-$\CC $ Column Partition in Example 2.}
		\label{fig:asbd-adbs}
	\end{minipage}
	\begin{minipage}{0.495\linewidth}
		\centering
		\includegraphics[width=1\linewidth]{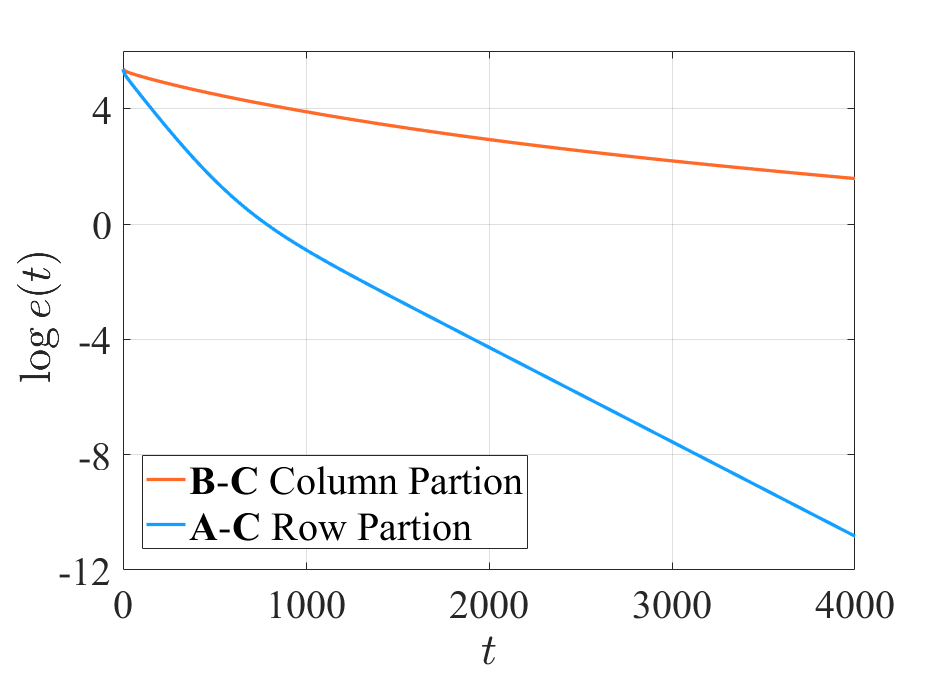}
		\caption{The trajectories of $ \log e(t) $ for  problem   $ \{\A^{[1]},\B^{[1]},\CC^{[1]}\} $ under partitions $ \B$-$\CC $ Column and $ \A$-$\CC $ Row, respectively, in  Example 2.}
		\label{fig:sd-ds-}
	\end{minipage}
\end{figure}

In Fig. \ref{fig:asbd-adbs}, we plot the trajectories of $ e(t) $ (as setting $ K=1 $ in \eqref{eq:et}) in logarithmic scales under the $ \B$-$\CC $ Column Partition for the two data sets $ \{\A^{[1]},\B^{[1]},\CC^{[1]}\} $ and $ \{\A^{[2]},\B^{[2]},\CC^{[2]}\} $, respectively. It can be seen that faster convergence is achieved at $ \{\A^{[2]},\B^{[2]},\CC^{[2]}\} $.
In  Fig. \ref{fig:sd-ds-}, we plot the trajectories of $ e(t) $ in logarithmic scales for $ \{\A^{[1]},\B^{[1]},\CC^{[1]}\} $ under partitions $ \B$-$\CC $ Column and $ \A$-$\CC $ Row, respectively.  These figures show that the size of data for each node in different partitions has an effect on the  convergence rate. Such examples motivate us to advance a conjecture about the existence of data complexity vs. convergence speed tradeoffs for the design of distributed algorithms.

\subsection{Discussions}
\subsubsection{General Sylvester Equation}\label{subsec:genSE}
Consider the Sylvester equation in its general form:
\begin{equation}\label{eq:ori_Sylvester}
\A\X+\X\B=\CC, \ \A\in\R^{n\times n}, \B\in\R^{m\times m}, \CC\in\R^{n\times m}.
\end{equation}	
Note that the vectorized form \eqref{eq:olinear_equation} continues to apply to \eqref{eq:ori_Sylvester}.  
We  define a $ m $-node ($ n $-node) network  under the $ \B$-$\CC $ Column ($ \A$-$\CC $ Row) Partition.  
Then the flow \eqref{eq:Flow-n-nodes} can be utilized in the same form, leading to the  convergence results under slightly different indices, e.g., under the $ \B$-$\CC $ Column Partition the limit of rate in Theorem \ref{thm:convergence} will be read as
$$ \lim\limits_{K\to \infty}r(K)= \lambda_{\rank(\Hh)}(\frac{1}{m}(\sum_{i=1}^{m}\Hh^{\dag}_{i}\Hh_{i})). $$
\subsubsection{Higher-resolution Data Partition}\label{subsubsec:n^2}
Define a $ n^2$-node network with node set $ \V^{\mathrm{H}}=\{1,2,\cdots,n^2\} $ forming a  graph $ \G^{\mathrm{H}}=(\V^{\mathrm{H}}, \EE^{\mathrm{H}}). $
Suppose that the index of node $ i $ satisfies $ i=(k-1)n+l $ with $ k=1,\cdots,n, l=1,\cdots,n. $ Then any $ i\in \V^{\mathrm{H}} $ can be uniquely represented by a binary array $ (k,l). $
Here we have the partition that node $ i(k,l) $ holds 
the $ l $-th row of $ \A $, the $ k$-th column of $ \B $ and 
the $ (l,k)$-th entry of $ \CC $, denoted by
$ (\A^T)_l^T, B_k $ and $ \CC_{lk} $.
Denote $$ [\A_0]_{i(k,l)}:=e_k^T\otimes(\A^T)_l^T,\quad [\B_e]_{i(k,l)}:=e_l^T(\B_k^T\otimes \I_n). $$ 
Then node $ i(k,l) $ has access to the equation
$ \big([\A_0]_{i(k,l)}+[\B_e]_{i(k,l)}\big)\mathrm{vec}(\X)=\CC_{lk}. $
For $ i(k,l)=i\in\V^{\mathrm{H}} $, denote
\begin{align*}
\E^{\mathrm{H}}_{i(k,l)}:=\{\y\in \R^{n^2}:\Big([\A_0]_{i(k,l)}+[\B_e]_{i(k,l)}\Big)\y=\CC_{lk}\}.
\end{align*} 
Case  (\uppercase\expandafter{\romannumeral1}) and Case  (\uppercase\expandafter{\romannumeral2}) guarantee that $ \E^{\mathrm{H}}_{i(k,l)} $ and their
intersection are nonempty. Therefore, our preceding discussion can be easily applied to this partition with designing the flow
\begin{equation*}
\dot{\x}_{i}=K(\sum_{j\in \N_i}(\x_j-\x_i))+\Pp_{\E^{\mathrm{H}}_{i(k,l)}}(\x_i)-\x_i, \quad i\in \V^{\mathrm{H}}.
\end{equation*}
With $ \E^{\mathrm{H}}:=\cap_{i=1}^{n}\E^{\mathrm{H}}_{i(k,l)}$ and $$ \h_i:=[\A_0]_{i(k,l)}+[\B_e]_{i(k,l)}\in \R^{1\times n^2}$$ $ (\h_i^{\dag}\h_i=\frac{\h_i^T\h_i}{\h_i\h_i^T} $ if $ \h_i\ne \0 , \h_i^{\dag}\h_i=\0 $ otherwise), we have
\begin{equation*}\label{eq:converge-xi-n^2-nodes}
\lim\limits_{t\to \infty}\x_i(t)=\sum_{i=1}^{n^2}\Pp_{\E^{\mathrm{H}}}(\x_i(0))/n^2, \quad \forall i\in \V^{\mathrm{H}}.
\end{equation*}
The rate of exponential convergence $ \tilde{r}(K) $ satisfies that  $$ \lim\limits_{K\to\infty}\tilde{r}(K)= \lambda_{\rank(\Hh)}(\frac{1}{n^2}(\sum_{i=1}^{n^2}\h^{\dag}_i\h_i)). $$ 
Compared with the case for an $ n $-node network, in which each node  holds $ n^2+2n $ scalar elements of data,  each node  only needs to hold $ 2n+1 $ scalar elements of data in the higher-resolution data partition case for an $ n^2 $-node network. 

\begin{figure}[h!]
	\begin{minipage}{0.495\linewidth}
		\centering
		\includegraphics[width=1\linewidth]{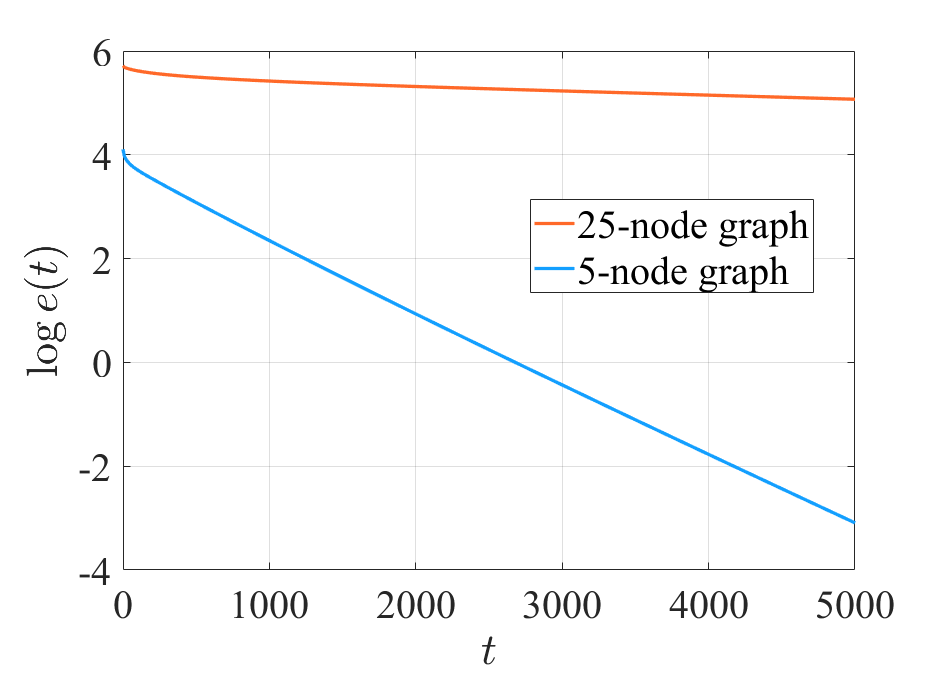}
		\caption{The trajectories of $ \log e(t) $ for two networks with cycle graph and $ K=1 $ in  Example 3.}
		\label{fig:cycle_vs}
	\end{minipage}
	\begin{minipage}{0.495\linewidth}
		\centering
		\includegraphics[width=1\linewidth]{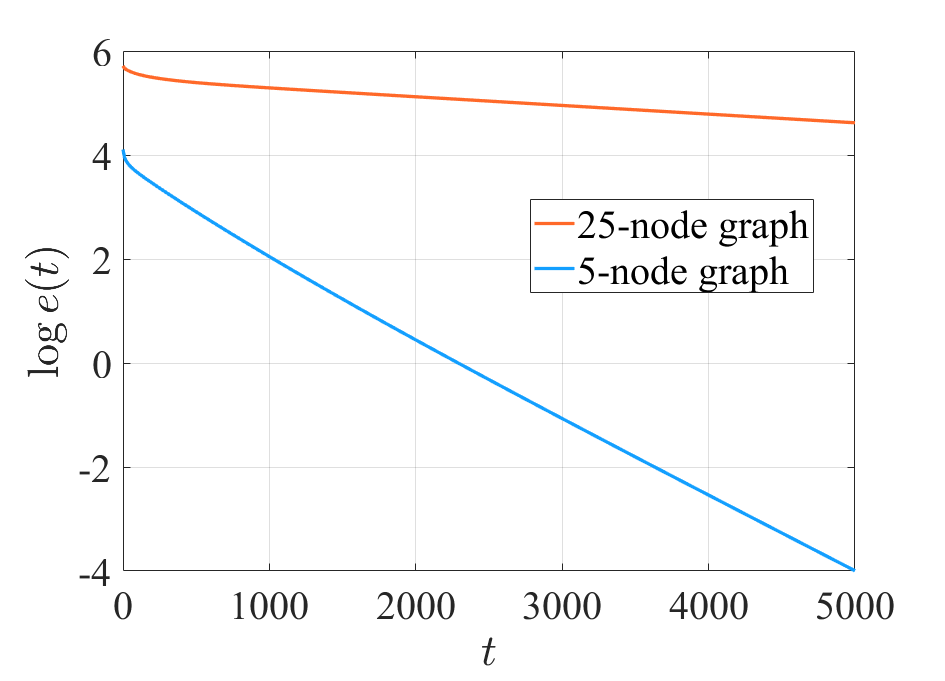}
		\caption{The trajectories of $ \log e(t) $ for two networks with complete graph and $ K=1 $ in Example 3.}
		\label{fig:comp_vs}
	\end{minipage}
\end{figure}
\noindent{\bf Example 3. }
Consider the same matrix equation as in {\rm Example 1}. Setting  $ K=1 $ in \eqref{eq:et}, we plot the trajectories of $ \log e(t) $ in Figs. \ref{fig:cycle_vs} and \ref{fig:comp_vs}, where we adopt  cyclic graphs and  complete graphs with $ 5 $ and $ 25 $ nodes, respectively. Each node in the  $ 5$-node network holds the data with the size of $ (5^2+2\times 5) $, while each node in the $ 25$-node network holds $ (2\times 5+1) $. Though the data complexity does decrease for every node in the $ 25$-node network, the convergence rate of the $ 5 $-node network is much faster.
 This indicates the existence of data distribution vs. convergence speed tradeoffs for the design of distributed algorithms. 

\subsubsection{Lyapunov Equations}
When $ \B=\A^T, $ we consider a Lyapunov equation with respect to variable $ \X\in \R^{n\times n} $:
\begin{equation}\label{eq:Lyap eq}
\A\X+\X\A^T=\CC,\quad \A,  \CC\in \R^{n\times n},
\end{equation}
where $ \CC $ is a symmetric matrix. If $ \X_0 $ is  the unique solution to \eqref{eq:Lyap eq}, it must hold that $ \X_0=\X_0^T. $ If $ \X_0 $ is a solution to \eqref{eq:Lyap eq}  when there exist an  infinite number of solutions, it must hold that $ \X_0^T $ is also a solution to \eqref{eq:Lyap eq}.
Due to the symmetry of $ \CC $, under the higher-resolution data partition in \ref{subsubsec:n^2}, we can alternatively adopt a network $ \hat{\G}=(\hat{\V},\hat{\mathrm{E}}) $ with $ n(n+1)/2 $ nodes rather than $ n^2 $ nodes. The node $ i $ is assigned with the data set
\begin{align*}
\mathcal{F}_i:=\big\{&  \big((\A^T)_l^T,  (\A^T)_k^T, \CC_{lk}, \CC_{kl}\big): k\le l, k, l\in\{1,\cdots,n\}i=g(k,l)=(k-1)n+l-k(k-1)/2 \big\}.
\end{align*} 
Defining $ f(k,l)=(k-1)n+l $, we introduce the affine subspaces 
\begin{align*}
\hat{\E}_{i}:=\Big\{&\y\in \R^{n^2}:\hat{\Hh}_i\y=\col\{\CC_{lk},
\CC_{kl}\}, \hat{\Hh}_i=\begin{bmatrix}
[\A_0]_{f(k,l)}+[\A^T_e]_{f(k,l)}\\
[\A_0]_{f(l,k)}+[\A^T_e]_{f(l,k)}
\end{bmatrix},\\& i=g(k,l), k\le l, k, l\in\{1,\cdots,n\}\Big\},\quad i\in\hat{\V},
\end{align*}
where $ [\A_0]_{f(k,l)}:=e_k^T\otimes(\A^T)_l^T $ and  
$[\A^T_e]_{f(k,l)}:=e_l^T(\A_k\otimes \I_n) $.
We can  modify the flow \eqref{eq:Flow-n-nodes} to
\begin{equation}\label{eq:Flow-LYn-nodes}
\dot{\x}_i=K(\sum_{j\in \N_i}(\x_j-\x_i))+\Pp_{\hat{\E}_i}(\x_i)-\x_i,\quad i\in\hat{\V}.
\end{equation}
Based on the same analysis, along \eqref{eq:Flow-LYn-nodes}  $ \vvec^{-1}(\x_{i}(t)) $ continues to converge to a solution of \eqref{eq:Lyap eq}, with the
rate of exponential convergence described by  $ \hat{r}(K) $, and $$ \lim\limits_{K\to\infty}\hat{r}(K)= \lambda_{\rank(\Hh)}(\frac{1}{n(n+1)/2}(\sum_{i=1}^{n(n+1)/2}\hat{\Hh}^{\dag}_i\hat{\Hh}_i)). $$

\section{Full Row/Column  Partition }\label{sec:row+column partition}

In this section, we investigate the full partitions of the data matrices $ \A $, $ \B $ and $ \CC $ along rows and columns, and present effective flows to solve the equation \eqref{eq:Sylvester} under such partitions over an $ n $-node network.
\subsection{Full Row/Column Partition }
We consider two full  partitions of the $ (\A,\B,\CC)$-triplet. 
\begin{itemize}
	\item[(\romannumeral1)] $[ \A $ Row/$ \B $-$ \CC $ Column Partition$ ] $ Node $ i $ holds the $ i $-th row of $ \A $, and the $ i $-th column of $ \B $ and $ \CC $, denoted by $ (\A^T)_i^T, \B_i $ and $ \CC_i $, respectively.
	\item[(\romannumeral2)] [$ \A $-$ \CC $  Row/$ \B $ Column Partition]  Node $ i $ holds the $ i $-th row of $ \A $ and $ \CC $, and the $ i $-th column of $ \B $, denoted by $ (\A^T)_i^T, (\CC^T)_i^T $ and $ \B_i $, respectively.
\end{itemize}
These two partitions are equivalent in the sense of $ \A\X+\X\B=\CC $ and $ \B^T\X^T+\X^T\A^T=\CC^T. $ Thus, we  focus on the $\A $ Row/$ \B $-$ \CC $ Column Partition in our analysis.
\begin{rem}
There may be other full row/column partitions, obtained e.g. through partitioning $ \A, \B $ by row and $ \CC $ by column or partitioning  $ \A , \B , \CC $ by column. By using appropriate equivalence transformation, we can deal with these partitions in a similar way, so more specific details are omitted.
 \end{rem} 

\subsection{An Augmented ``Consensus $ + $ Projection'' Flow}
For $ i\in\V $, define 
\begin{equation}\label{eq:spaceEz_i}
\begin{aligned}
&\E_i^{\text{Aug}}=\{\y=\vvec([\X,\mathbf{Z}])\in\R^{n^2(1+n)}: \ee_i(\A^T)_i^T\X+\X\B_i\ee_i^T-(({\LL_{\G}^T})_i^T\otimes\I_{n})\mathbf{Z}=\CC_{i}\ee_i^T \},
\end{aligned} \end{equation} 
 and a projection mapping $ \Pp_{\E_i^{\text{Aug}}}: \R^{n^2(1+n)}\to\E_i^{\text{Aug}} $.
We propose the following augmented network flow for $ \y_i(t)=\col\{\x_i(t),\z_i(t)\} $ with $ \x_i(t)\in\R^{n^2} $, 
\begin{equation}\label{eq:row+Flow-n-nodes}
\dot{\y}_i
=K\Big(\sum_{j\in \N_i}(\y_j-\y_i)\Big)+\Pp_{\E_i^{\text{Aug}}}(\y_i)-\y_i, \quad i\in\V,
\end{equation}
where $ \y_i\in\R^{n^2(1+n)} $, and $ \vvec^{-1}(\x_i)\in\R^{n\times n} $ is what we are interested in.
\begin{thm}\label{thm:row+column convergence}
	Under the  $\A $ Row/$ \B $-$ \CC $ Column Partition, for any initial value $ \y_0=\col\{\y_1(0),\cdots,\y_n(0)\} $,  there exists $ \X^*(\y_0)\in\R^{n\times n} $ as a solution to \eqref{eq:Sylvester}, such that  $ \vvec^{-1}(\x_i(t))=\vvec^{-1}([\I_{n^2},\0_{n^2\times n^3}]\y_i(t)) $ along the flow \eqref{eq:row+Flow-n-nodes} converges to $ \X^*(\y_0) $ exponentially, for all $ i\in\V $. 
\end{thm}
The proof of Theorem \ref{thm:row+column convergence} is given in Appendix  \ref{App:full}.

\subsection{Application for the Motivating Example}
The network flow \eqref{eq:row+Flow-n-nodes} can be used to solve the problem arising from the motivating example mentioned in subsection \ref{subsec:motivating-ex}.
Each node $ i $ representing subsystem $ i $ only knows the information $ \D_{ij}\in\R^{m\times m}, j\in\V $ and communicates with its neighbors for exchanging state information.
We utilize the flow \eqref{eq:row+Flow-n-nodes} via substituting 
\begin{align*}
&(\A^T)_i^T=[\D_{i1},\cdots\D_{in}],\ \B_i=\col\{\D_{i1},\cdots\D_{in}\},\  \CC_i=-(\I_{mn})_i,
\end{align*}  into  $ \E_i^{\text{Aug}} $ in \eqref{eq:spaceEz_i}.
As a result, along the flow \eqref{eq:row+Flow-n-nodes} carried out over $ \G $,  node $ i $ can obtain an evolutionary state $ \vvec^{-1}(\x_i(t)) $ by communicating and computing. Then every node can draw a conclusion about the stability of the overall network after  confirming two conditions:
\begin{itemize}
	\item[(\romannumeral1)] Each $ \vvec^{-1}(\x_i(t)) $ converges to a positive definite matrix at node $ i $;
	\item[(\romannumeral2)] All $ \vvec^{-1}(\x_i(t)) $ converge to the same limit.
\end{itemize}
It is easy to see while condition (\romannumeral1) can be verified by each node $ i $ by itself, and condition (\romannumeral2) can be established distributedly by for example, running a consensus algorithm for the node state limits.

\begin{figure}[h!]	
	\centering
	\includegraphics[width=0.3\linewidth]{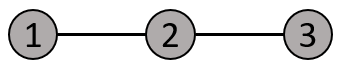}
	\caption{\small{The communicate graph in Example 4.}}
	\label{fig:3-nodes}
\end{figure}

\begin{figure}[htbp]
	\centering
	\includegraphics[width=0.5\linewidth]{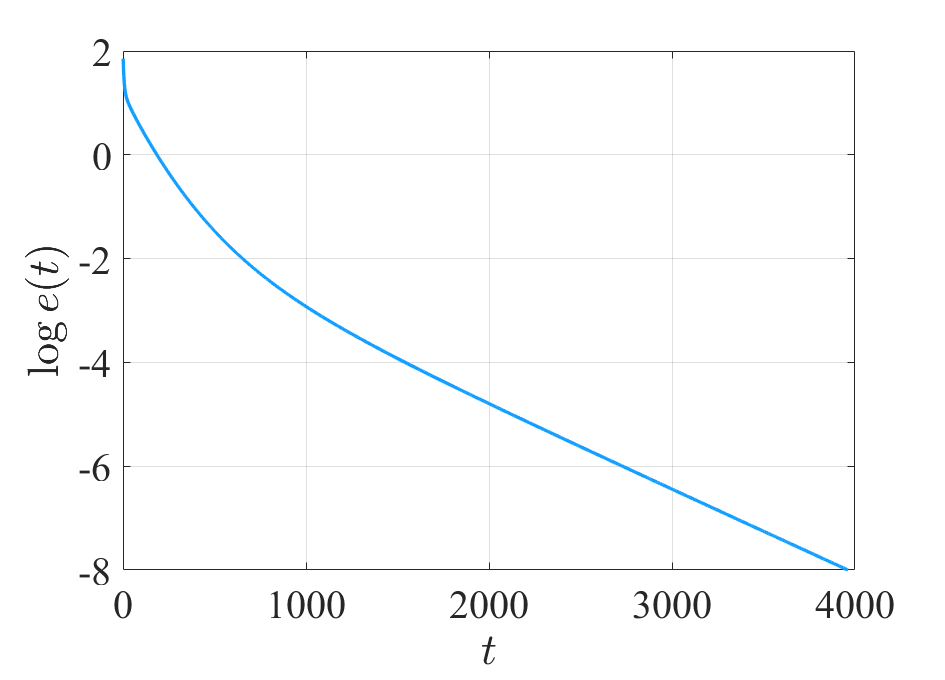}
	\caption{\small{The trajectories of $ e(t) $ in Example 4.} }
	\label{fig:fullloget}
\end{figure}

\noindent{\bf Example 4. } Consider three subsystems $  i=1,2,3 $, in 
the dynamics is in the form of \eqref{eq:subsys} with
\begin{align*}
\D_{11}&=\begin{bmatrix}
-5& 0\\0&-5
\end{bmatrix},\ \D_{12}=\begin{bmatrix}
4&1\\0&4
\end{bmatrix},\ \D_{13}=\0_{2\times 2},\\
 \D_{21}&=\begin{bmatrix}
1&0\\0&1
\end{bmatrix},\ \D_{22}=\begin{bmatrix}
-6& 1\\1&-3
\end{bmatrix},\ \D_{23}=\begin{bmatrix}
2&0\\2&4
\end{bmatrix},\\
 \D_{31}&=\0_{2\times 2},\ \D_{32}=\begin{bmatrix}
2& 0\\1&-1
\end{bmatrix},\ \D_{33}=\begin{bmatrix}
-4& 0\\0&-4
\end{bmatrix}.
\end{align*} 
The network communication structure is shown as Fig. \ref{fig:3-nodes}, whose Laplacian matrix is 
\[
\LL_{\G}=[1, -1, 0; -1,2,-1;0,-1,1].\]
Each node $ i $ can compute 
$  \Pp_{\E_i^{\text{Aug}}}(\cdot) $ in \eqref{eq:spaceEz_i} according to Lemma \ref{lemma:projection}. Under the network flow \eqref{eq:row+Flow-n-nodes}, there holds
\begin{align*}\small
&\lim\limits_{t\to\infty}\vvec^{-1}(\x_i(t))=\mathbf{P}^*=\begin{bmatrix}
0.2278&	 0.1343& 	0.1176&	0.1690&	0.0744&	-0.0009\\
0.1343&	0.3170&	0.0990&	0.2713&	0.0694&	-0.0068\\
0.1176&	0.0990&	0.1529&	0.1360&	0.0819&	0.0040\\
0.1690&	0.2713&	0.1360&	0.4106&	0.1067&	0.0278\\
0.0744&	0.0694&	0.0819&	0.1069&	0.1660&	-0.0021\\
-0.0009&	-0.0068&	0.0040&	0.0278&	-0.0021&	0.1190
\end{bmatrix},
\end{align*}
where $ \mathbf{P}^* $ is a positive definite solution to $$ \A\mathbf{X}+\mathbf{X}\A^T=-\I_{6} $$ with $ \A=[\D_{11},\D_{12},\D_{13};\D_{21},\D_{22},\D_{23};\D_{31},\D_{32},\D_{33}].$ See from Fig. \ref{fig:fullloget} the trajectory of $$ e(t)=\frac{1}{3}\sum_{i=1}^3\|\vvec^{-1}(\x_i(t))-\mathbf{P}^*\|_F^2 $$ in logarithmic scales.
 As for every subsystem, on the one hand, it can hold the information that $ \vvec^{-1}(\x_i(t)) $ converges to a positive definite matrix $ \mathbf{P}^* $. On the other hand, they can carry out a consensus test and have confirmed that all the subsystems states along \eqref{eq:row+Flow-n-nodes} are reaching a consensus state. Then every subsystem can conclude that the whole system is stable, which is in agreement with the fact that the global matrix $ \A $ is Hurwitz.

\section{Clustering Block Partition}\label{sec:double partition}
In this section, we turn to clustering block partitions of $ \A, \B, \CC $ for seeking distributed solutions of the equation \eqref{eq:Sylvester}. 
It seems  possible that, for certain structured matrices, general block partitions may be particularly useful. 
\subsection{Clustering Block Partition}\label{subsec:double}
Consider
a double-layer network that has $ n $ clusters with each cluster having $ n $ nodes. These clusters are indexed in $ \V=\{1,\cdots,n\} $ forming an outer layer graph $ \G=(\V,\EE), $ while the nodes in cluster $ i $ are indexed in $ \V_i=\{i_1,\cdots,i_n\} $ forming an inner layer graph $ \G_i=(\V_i,\EE_i). $ In total there are $ n^2 $ nodes in the overall network. The neighbor set of cluster $ i $ is given by $ \N_i:=\{j: (i,j)\in \EE\} $, which means that nodes in cluster $ i $ can receive  information from  nodes in its neighbor clusters; meanwhile, the neighbor set of node $ i_j $ in cluster $ i $ is given by $ \N_{i_j}:=\{i_k: (i_j,i_k)\in \EE_i\} $, which means that node $ i_j $ can receive  information from its neighbor nodes in its own cluster. Let $ \LL_{\G} $ and $ \LL_{\G_i} $ denote the Laplacian matrix of the outer layer graph (linking the clusters) and inner layer graphs (linking nodes in each cluster), respectively.
We recall $ \Hh=\I_n\otimes \A+\B^T\otimes \I_n $ and define the $ i $-th column of $ \CC $ as 
$ \CC_i=\sum\nolimits_{j=1}^n\CC_{ji}\ee_j $ with $ \CC_{ji} $ being the $ (j,i) $-th entry of $ \CC $.
 Define an indicator function $ \1_{\{j=i\}} $, where 
$ \1_{\{j=i\}}=1 $ if $ j=i, $ and $ \1_{\{j=i\}}=0, $ otherwise.
Then we  consider the following data partition. 

\noindent [Column  $ \B$-$\CC $ Block Partition, as in Table \ref{tab:double-data}]  
The node $ i_j $ holds $ \B_{ji} $ and $ \CC_{ji} $ (the $ (j,i) $-th entry of $ \B $ and $ \CC $), and additionally, the node $ i_i $ holds $ \A $. Together, we say
cluster $ i $ holds $ \A, $ $ \B_i $ and $ \CC_i $.
Specifically, each  node $ i_j $ holds a state $ \x_{i_j}\in \R^n $, while the cluster state $ \x_i=\col\{\x_{i_1},\cdots,\x_{i_n}\}\in \R^{n^2} $ satisfies
$$ \sum\limits_{j=1}^{n}\left((\1_{\{j=i\}}\A+\B_{ji}\I_n)\x_{i_j}-\CC_{ji}\ee_j\right)=\0_n. $$
\begin{table}[h!]
	\centering	\caption{The Column  $ \B$-$\CC $ Block Partition.}
	\label{tab:double-data}
\begin{tabular}{|c|c|c|c|c|}
	\hline 
	& $ 1 $-st Node & $ 2 $-nd Node & $\cdots $ & $ n $-th Node \\ 
	\hline 
	Cluster 1 & $ \A, \B_{11}, \CC_{11} $ & $ \B_{21}, \CC_{21} $ & $\cdots $ & $  \B_{n1}, \CC_{n1} $ \\ 
\hline 
Cluster 2 & $ \B_{12}, \CC_{12} $ & $ \A, \B_{22}, \CC_{22} $ & $\cdots $ & $  \B_{n2}, \CC_{n2} $ \\ 
\hline 
$\vdots $ & $\vdots $ & $\vdots $ & $\ddots $ & $\vdots $ \\ 
\hline 
Cluster n & $  \B_{1n}, \CC_{1n} $ & $ \B_{2n}, \CC_{2n} $ & $\cdots $ & $ \A, \B_{nn}, \CC_{nn} $ \\ 
\hline 
\end{tabular} 
\end{table}

Therefore, all the estimates from clusters need to reach a consensus $ \x_1^*=\cdots=\x_n^* $, which is the estimation of a solution to \eqref{eq:olinear_equation}. In view of \cite{wang2019scalable}, essentially any data block partition would work if the algorithm can be correspondingly designed.

\subsection{``Local Conservation + Global Consensus'' Flow}\label{subsec:double flow}
Take an auxiliary variable $ \z_{i_j}\in \R^n $ associated with and known by node $ i_j. $ Each node $ i_j $ can obtain the information about $ \z_{i_k}, i_k\in \N_{i_j} $ from its neighbors within the same cluster. Then the auxiliary variables of all nodes within the cluster $ i $ combine the cluster variable $ \z_{i}=\col\{\z_{i_1},\cdots,\z_{i_n}\}\in \R^{n^2} $. Each node $ i_j $ holds state $ \x_{i_j}\in\R^n $ and gets the information about $ \x_{k_j}, k\in \N_{i} $ from its neighbor clusters and then the states of all nodes within the cluster $ i $ combine the cluster state $ \x_{i}=\col\{\x_{i_1},\cdots,\x_{i_n}\}\in \R^{n^2} $. Let $ K>0 $ be a given constant. We propose the following continuous-time network flow:
\begin{equation}\label{eq:double-flow 1st partition_xi_j}
\begin{aligned}
\dot{\x}_{i_j}=&-(\1_{\{j=i\}}\A+\B_{ji}\I_n)^T((\1_{\{j=i\}}\A+\B_{ji}\I_n)\x_{i_j}-\CC_{ji}\ee_j-\sum_{i_k\in \N_{i_j}}(\z_{i_j}-\z_{i_k}))-K\sum_{k\in \N_{i}}(\x_{i_j}-\x_{k_j}),\\
\dot{\z}_{i_j}=&(\1_{\{j=i\}}\A+\B_{ji}\I_n)\x_{i_j}-\CC_{ji}\ee_j-\sum_{i_k\in \N_{i_j}}(\z_{i_j}-\z_{i_k}).
\end{aligned}
\end{equation}
Denote $ \M_i=\diag\{\1_{\{j=i\}}\A+\B_{ji}\I_n, j=1,\cdots,n\} $ and $ \tilde{\CC}_i=\col\{\CC_{1i}\ee_1,\cdots,\CC_{ni}\ee_n\}$; then we reformulate \eqref{eq:double-flow 1st partition_xi_j} as
\begin{equation}\label{eq:double-flow 1st partition_x_i}
\begin{aligned}
\dot{\x}_{i}=&-\M_i^T(\M_i\x_{i}-\tilde{\CC}_{i}-(\LL_{\G_i}\otimes\I_n)\z_i)-K\sum_{k\in \N_{i}}(\x_{i}-\x_{k}),\\
\dot{\z}_{i}=&\M_i\x_{i}-\tilde{\CC}_{i}-(\LL_{\G_i}\otimes\I_n)\z_i.
\end{aligned}
\end{equation}
We further define $ \bar{\M}=\diag\{\M_i, i=1,\cdots,n\} $, $ \bar{\CC}=\col\{\tilde{\CC}_{i},\cdots,\tilde{\CC}_{n}\} $, $ \bar{\LL}=\diag\{\LL_{\G_i}\otimes\I_n,i=1,\cdots,n\} $ and $ \z=\col\{\z_1,\cdots,\z_n\}, \x=\col\{\x_1,\cdots,\x_n\} $.  The flow \eqref{eq:double-flow 1st partition_x_i} can be rewritten  as a compact form
\begin{equation}\label{eq:double-flow 1st partition_x}
\begin{aligned}
\dot{\x}&=-\bar{\M}^T(\bar{\M}\x-\bar{\CC}-\bar{\LL}\z)-K(\LL_{\G}\otimes\I_{n^2})\x,\\
\dot{\z}&=\bar{\M}\x-\bar{\CC}-\bar{\LL}\z.
\end{aligned}
\end{equation}
Denoting $ \Gg:=\begin{bmatrix}
\bar{\M}^T\bar{\M}+K(\LL_{\G}\otimes\I_{n^2}) & -\bar{\M}^T\bar{\LL}\\
-\bar{\M} & \bar{\LL}
\end{bmatrix} $,
 we present the following theorem.

\begin{thm}\label{thm:doublecase}
	Under the Clustering Block Partition, for any initial values $ \x_0=\col\{\x_1(0),\cdots,\x_n(0)\},$ and $ \z_0=\col\{\z_1(0),\cdots,\z_n(0)\} $, there exists $ \X^*(\x_0,\z_0)\in\R^{n\times n} $ as a solution to \eqref{eq:Sylvester}, such that  $ \vvec^{-1}(\x_i(t)) $ along the flow \eqref{eq:double-flow 1st partition_xi_j} converges to $ \X^*(\x_0,\z_0) .$ 	Moreover, there exist $ \beta(\x_0, \z_0), r^*(K)>0 $, such that for all $  t\ge0, $
	$$  \sum_{i=1}^n\|\vvec^{-1}(\x_i(t))-\X^*(\x_0,\z_0)\|_F^2\le \beta(\x_0,\z_0) e^{-2r^*(K)t}, $$
	 where the rate of exponential convergence $ r^*(K) $ is a non-decreasing function with respect to $ K $ satisfying $ r^*(K)=\lambda_{k}(\Gg) $, where $ k=\rank(\Gg)\le2n^3-n^2-\dim(\cap_{i=1}^n\ker(\M_i)) $.
\end{thm}
The proof of Theorem \ref{thm:doublecase} is shown in Appendix \ref{App:double-con}.

\subsection{Numerical Example}	
\begin{figure}[htbp]
		\centering
		\includegraphics[width=0.5\linewidth]{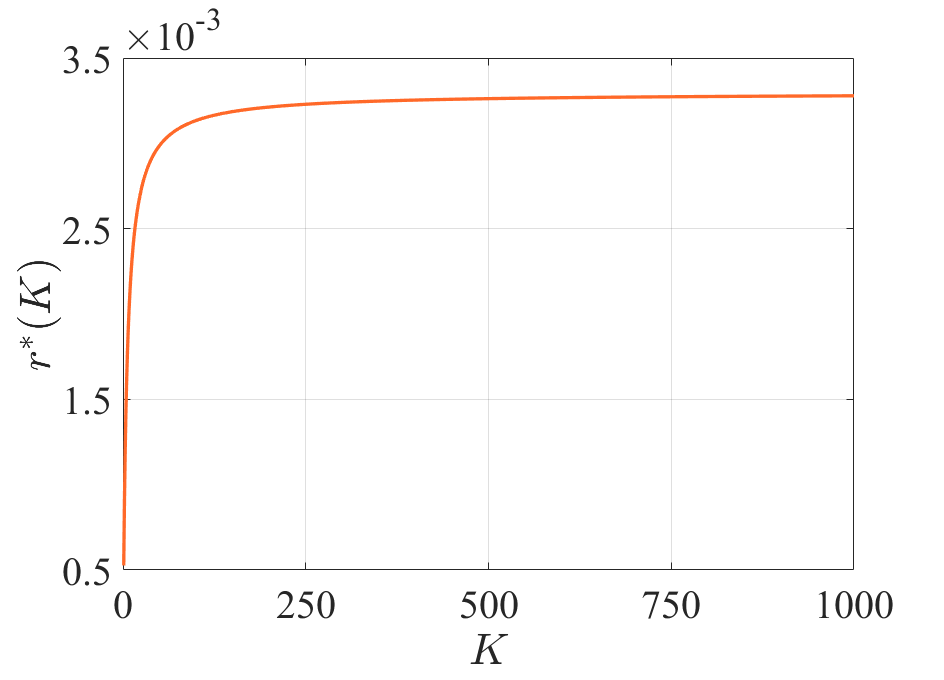}
		\caption{The trajectory of  $ r^*(K) $ over the parameter $ K $.}
		\label{fig:rstark}
\end{figure}

\noindent{\bf Example 5. } Consider the same matrix equation as in {\rm Example 1}, which has a unique solution. We use the two kinds of networks with a 25-node graph in  subsection \ref{subsubsec:n^2} and a graph of five 5-node clusters in  subsection \ref{subsec:double}, respectively. 
Define the error function under the Column  $ \B$-$\CC $ Block Partition:
\begin{align*}
E_{K}^d(t)&=\sum_{i=1}^{5}\|\x_{i}(t)-\mathrm{vec}(\X^*)\|^2=\sum_{i=1}^{5}\|\col\{\x_{i_1},\cdots,\x_{i_5}\}-\mathrm{vec}(\X^*)\|^2.
\end{align*}  	
\begin{figure}[h!]
	\centering
	\includegraphics[width=0.5\linewidth]{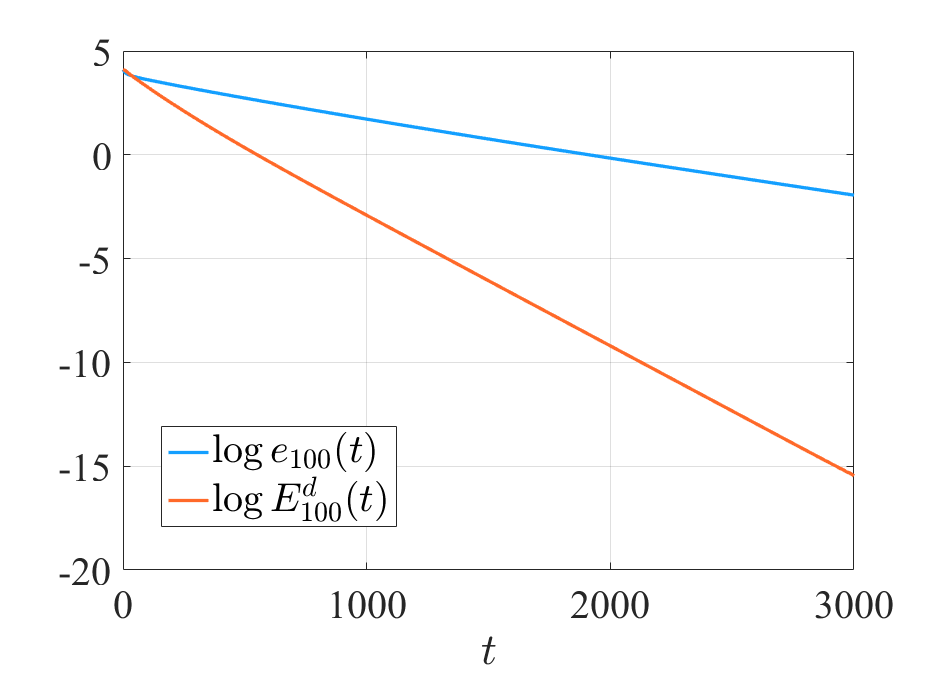}
	\caption{The trajectories of $ \log E_{100}^d(t) $ and $ e_{100}(t) $, respectively.}
	\label{fig:doublevs}
\end{figure}
For a complete graph,
taking the zero matrix as the initial value,  we plot in Fig. \ref{fig:doublevs} the trajectories of $ E_{100}^d(t) $ and $ e_{100}(t) $ defined in \eqref{eq:et} in logarithmic scales for $ \x_{i_j}(t) $ evolving along \eqref{eq:double-flow 1st partition_xi_j} and $ \x_i(t) $  along \eqref{eq:Flow-n-nodes}, respectively. Fig. \ref{fig:doublevs} shows that $ \x_{i_j}(t) $ along the flow \eqref{eq:double-flow 1st partition_xi_j} converges exponentially and the convergence rate of clustering block partition is much faster than that of  partial  $ \B$-$\CC $ Column Partition in Section \ref{sec:row/column partition}.
With different values of  $ K, $ we can also calculate $ r^*(K) $ and plot the trajectory of $ r^*(K) $ over $ K $ in Fig. \ref{fig:rstark}, which shows that the rate of exponential convergence is a non-decreasing function with respect to $ K $.
Results in these figures are consistent with Theorem \ref{thm:doublecase}.

\section{Conclusion}\label{sec:conclusion}
This paper has focused on the distributed computation of the multi-agent network for Sylvester matrix equations.  We have proposed  several network flows for partitions of partial  row/column, full row/column  and clustering block  about the data matrices,  inspired by the computation for linear algebraic equations. We have remarked on a special case for symmetric solutions and discussed the general Sylvester equation and others with examples. Accordingly, appropriate partitions could be selected based on actual conditions. The convergence and the limitation of convergence rate have been established in view of matrix theory and linear system theory, which also have been verified by typical numerical examples.    Future work includes characterizing the performance of the proposed solutions for general directed and switching networks, as well as methods for accelerating the flows by optimizing network structures.

\medskip

\medskip

\noindent{\bf \Large Appendices}

\appendix
\section{Preliminaries}\label{sec:preliminaries}
In this appendix, we present some preliminaries on matrix analysis, affine spaces, and exponential stability of dynamical systems.

	For a matrix $ \A \in \R^{n\times m}, $ a M-P pseudoinverse \cite{penrose1955generalized} of $ \A $ is defined as a matrix $ \A^{\dag}\in \R^{m\times n} $
	satisfying all of the following four equalities:
	$ (\romannumeral1)\A\A^{\dag}\A=\A; (\romannumeral2)\A^{\dag}\A\A^{\dag}=\A^{\dag};
	(\romannumeral3)(\A\A^{\dag})^T=\A^{\dag}\A;
	(\romannumeral4)(\A^{\dag}\A)^T=\A^{\dag}\A. $ Then the following lemma about the pseudoinverse holds, as well as a lemma about the inequalities of eigenvalues.

\begin{lem}[\cite{bernstein2005matrix}]\label{lemma:prop of inverse}
	For any matrix $ \A \in \R^{n\times m}, $ the following statements hold.
	\begin{enumerate}
		\item[(\romannumeral1)] The M-P pseudoinverse of matrix $ \A $ is unique. The pseudoinverse of the pseudoinverse is the original matrix: $ (\A^{\dag})^{\dag}=\A; $
		\item[(\romannumeral2)] $ \rank(\A^{\dag})=\rank(\A)=\rank(\A^T)=\rank(\A^{\dag}\A)=\rank(\A\A^{\dag}); $
		\item[(\romannumeral3)] $ \ker(\A^{\dag})=\ker(\A^T), $ $ \ker(\A)= \ker(\A^{\dag}\A); $
		\item[(\romannumeral4)] $ (\A^{\dag}\A)^2=\A^{\dag}\A=(\A^{\dag}\A)^T, $ $ \A^{\dag}\A $ is real symmetric and idempotent, and its eigenvalues can only be zero or one.
	\end{enumerate}
\end{lem}
\begin{lem}[Weyl's inequality \cite{franklin2012matrix}]\label{lemma:Weyl's inequality}
	Suppose that $ \M $ and $ \NN $ are $ n\times n $ symmetric matrices. Then
	\begin{equation}\label{eq:Weyl's inequality}
	\lambda_j(\M)+\lambda_n(\NN)\le \lambda_{j}(\M+\NN)\le \lambda_j(\M)+\lambda_1(\NN).
	\end{equation}
\end{lem}

Next, for a non-defective matrix, the following lemma holds, where
	a matrix $ \M\in\R^{n\times n} $ is non-defective if it is diagonalizable. 

\begin{lem}[\cite{juang1989eigenvalue}]\label{lemma:Eigenvalue Derivatives}
	Let $ \M $ be a non-defective matrix depending on a parameter $ \rho. $ Suppose that the eigenvalue $ \lambda_{1} $ has a multiplicity $ k_1 $ ($ \lambda_{i}=\lambda_{1} $ for $ i=1,\cdots,k_1 $). Let 
	$ 	\X_1=[\x_{1},\cdots,\x_{k_1}]$ and $ \Y_1=[\y_{1},\cdots,\y_{k_1}] $
	represent the base vectors of the left and the right eigenvector space associated with the eigenvalue $ \lambda_{1} $ for $ \M(\rho_0) $, respectively, where the chosen bases satisfy $ \X_1^T\Y_1=\I_{k_1} $. Then, 
	for $ i=1,\cdots,k_1,  \lambda_{i}^{'}=\frac{\partial \lambda_{i}(\rho)}{\partial \rho}|_{\rho=\rho_0} $, there holds
	\begin{equation}\label{eq: derivative of eigenvalue}
	\big(\X_1^T\frac{\partial \M(\rho)}{\partial \rho}|_{\rho=\rho_0}\Y_1\big)\z=\lambda_{i}^{'}\z,
	\end{equation}
	where $ \M\mathbf{\varphi}=\lambda_{i}\mathbf{\varphi} $ and $ \mathbf{\varphi}=\Y_1\z $ with some $ \z\in\R^{k_1\times 1}. $ Equivalently, the eigenvalue derivatives $ \frac{\partial\lambda_{i}}{\partial \rho}, i=1,\cdots,k_1 $ are the eigenvalues of  matrix $ \X_1^T\frac{\partial \M(\rho)}{\partial \rho}|_{\rho=\rho_0}\Y_1. $
\end{lem}
\begin{lem}[Lemma 1 in \cite{wang2019scalable}]\label{lemma:G_eigenvalue>0or=0}
	Let	\[\Q=\begin{bmatrix}
	\Q_1^T\Q_1+\Q_2 & -\Q_1^T\Q_3\\
	-\Q_1 & \Q_3
	\end{bmatrix}, \]
	where all submatrices in $ \Q $ are real matrices, and $ \Q_2, \Q_3 $ are positive semi-definite. Then all eigenvalues of $ \Q $ are greater than or equal to 0. Moreover, if $ \Q $ has a zero eigenvalue, the zero eigenvalue must be non-defective. 
\end{lem}

An affine space $ \cite{rockafellar2015convex} $ is a set $ \Aa $ if $ (1-\theta)\x+\theta\y\in\Aa $ for any $ \x, \y\in\Aa $ and $ \theta\in\R $. A projection mapping on an affine subspace is a linear transformation, which assigns each $ \x\in\Aa $ to the unique element $ \Pp(\x)\in\Aa $ such that $ \|\x-\Pp_{\Aa}(\x)\|=\min_{\y\in\Aa}\|\x-\y\|. $ For the affine space and projection, the following lemma holds.

\begin{lem}[\cite{bernstein2005matrix}]\label{lemma:projection}
	Let $ \K:=\{\y\in \R^m: \Gg\y=\z, \Gg\in\R^{n\times m}, \z\in\R^{m}\} $ be an affine subspace. Denote $ \Pp_{\K}: \R^{m}\to \K $ as the projector onto $ \K $. Then $ \Pp_{\K}(\x)=(\I_m-\Gg^{\dag}\Gg)\x+\Gg^{\dag}\z $ for all $ \x\in \R^{m} $, where $ \Gg^{\dag}\in \R^{m\times n} $ is the M-P  pseudoinverse of $ \Gg $.
\end{lem}

Finally, we introduce a concept about 
exponential convergence.
	A solution of the system \[\dot{\x}(t)=f(\x,t), \quad \x\in\R^n,\ t\ge0,\] is termed to be exponentially convergent to $ B_{\delta} $ at rate $ r $ \cite{corless1990guaranteed} if there exists $ r>0 $, and for any initial condition $ \x_0\in\R^n, $ there exists $ c(\x_0)>0 $, such that for any solution $ \x(t) $ with $ \x(t_0)=\x_0\in \mathcal{D}, $ where $ \mathcal{D} $ is an open set containing the origin, there holds
	\begin{equation}\label{def:exp rate to B}
	\|\x(t)\|\le\delta+c(\x_0)e^{-r(t-t_0)}, \quad \forall t\ge t_0.
	\end{equation} If in addition, $ \delta=0 $, this solution is exponentially convergent to zero.

\section{Proof of Theorem \ref{thm:convergence}}\label{App:1con}
\subsection{Preliminary Lemmas}
Recall the following lemma on the flow \eqref{eq:Flow-n-nodes} from \cite{shi2017networkTAC}.
\begin{lem}[\cite{shi2017networkTAC}]\label{lemma:invariant}
	Let {\rm Assumption \ref{assumption-exact solution}} hold. Assume that $ y^* $ is an exact solution to \eqref{eq:olinear_equation} and $ r>0 $ is arbitrary. Define $ \mathcal{M}^*(r):=\{\y\in\R^{n^2}: \|\y-\y^*\|\le r\}. $ Then $ (\mathcal{M}^*(r))^n=\mathcal{M}^*(r)\times\cdots\times\mathcal{M}^*(r) $ is a positively invariant set along the flow \eqref{eq:Flow-n-nodes}.
\end{lem}

Next, we introduce the  notations $ O $ and $ \Theta $. For two functions $ g, h $ with $ h(\cdot)>0 $, denote
\begin{itemize}
	\item $ g(t) = O(h(t)) $ as $ t\to\infty $ if there exist $ c, t_0>0 $, such that $ |g(t)|\le ch(t) $ for all $ t\ge t_0 $;
	\item $ g(t) = \Theta(h(t)) $ as $ t\to\infty $ if there exist $ c_1, c_2>0 $ and $ t_0>0 $, such that $ c_1h(t)\le |g(t)|\le c_2h(t) $ for all $ t\ge t_0 $.	
\end{itemize}
Then the following lemma is based on some basic convergence properties of linear time-invariant systems.
\begin{lem}\label{lemma:LTIorder}
	Consider a
	linear time-invariant system 
	\begin{equation}\label{eq:lemma8lti}
	\dot{\x}(t)=-\mathbf{F}\x(t)+\alpha(t),\quad t>0,\ \x(t)\in\R^m,
	\end{equation}
	where $ \mathbf{F}\in \R^{m\times m} $ is positive semidefinite and $ \rank(\mathbf{F})=k\le m $. Suppose that $ \|\alpha(t)\|=O(e^{-rt}) $ as $ t\to\infty $ with  $ r>\min\{\lambda\in \mathrm{spec}(\mathbf{F}):\lambda\neq0 \}=\lambda_{*}. $  Then, for an initial condition $ \x_0 $, the following statements hold.
	\begin{enumerate}
		\item[(\romannumeral1)]  There exists a unique $ \z^*(\x_0)\in \ker(\mathbf{F}) $, such that $ \lim_{t\to\infty}\x(t)=\z^*(\x_0) $.
		\item[(\romannumeral2)] 
		For almost all initial conditions, $ \|\x(t)-\z^*(\x_0)\|=\Theta(e^{-\lambda_{*}t}) $.
	\end{enumerate}	
\end{lem}
The result of Lemma \ref{lemma:LTIorder} is trivial to establish when $ \mathbf{F} $ is 1 by 1. For $ m>1 $, using an orthogonal matrix $ \T $ for which $ \T^{T}\mathbf{F}\T $ is diagonal can be helpful to finish the proof. The details of proof are omitted for space limitations. 

Next, we establish a lemma on the convergence of $ \x_{ave}=\frac{1}{n}\sum_{i=1}^n\x_i $ along the flow \eqref{eq:Flow-n-nodes}.
\begin{lem}\label{lemma:step3part}
	Along the flow \eqref{eq:Flow-X}, $ \x(t) $ is the solution for given $ \x_0, $ and $ \bar{\x}(t)=\1_n\otimes\x_{ave}(t) $. Then, for any $ \delta>0, $ any $ t_0>0, $ there exists $ K_{\delta,t_0} $, such that 
	\begin{equation}
	\|\x(t)-\bar{\x}(t)\|\le\delta,\quad \forall  K>K_{\delta,t_0},\ t>t_0.
	\end{equation}
\end{lem}
\noindent$ Proof. $
Following from  Lemma \ref{lemma:invariant} and \cite{shi2017networkTAC}, for given $ \x_0, \|\x(t)-\1_n\otimes\y^*(\x_0)\| $ is always bounded; moreover, $ \|\x_i(t)\| $ is bounded for all node $ i $ and $ \|\x(t)-\bar{\x}(t)\| $ is bounded as well. 
According to Lemma \ref{lemma:projection},  $$  \Pp_{\E_i}(\x_{i})=(\I_{n^2}-\Hh^{\dag}_i\Hh_i)\x_{i}+\Hh^{\dag}_i\CC_i, $$ it can be easily calculated that 
\begin{align*}
\dot{\bar{\x}}(t)=&\1_n\otimes(\frac{1}{n}\sum_{i=1}^{n}(\Pp_{\E_i}(\x_i(t))-\x_i(t)))
=\1_n\otimes(\frac{1}{n}\sum_{i=1}^{n}(-\Hh^{\dag}_i\Hh_i\x_{i}(t)+\Hh^{\dag}_i\CC_i)),
\end{align*}
which leads to that $ \dot{\bar{\x}}(t) $ is bounded. Next, we consider the property of $ \|\x(t)-\bar{\x}(t)\|^2, $
\begin{align*}
&\frac{d}{dt}\|\x(t)-\bar{\x}(t)\|^2
=2\langle \x(t)-\bar{\x}(t),\dot{\x}(t)-\dot{\bar{\x}}(t) \rangle\\
=&2\langle \x(t)-\bar{\x}(t), -(K\LL_{\G}\otimes \I_{n^2}+\J)\x(t)+\Q_{\CC}-\dot{\bar{\x}}(t)\rangle\\
=&2\langle \x(t)-\bar{\x}(t), -(K\LL_{\G}\otimes \I_{n^2})(\x(t)-\bar{\x}(t))\rangle+\phi(t)\\
\le&-2K\lambda_{n-1}(\LL_{\G})\|\x(t)-\bar{\x}(t)\|^2+\phi(t),
\end{align*}
where $ \phi(t)=2\langle\x(t)-\bar{\x}(t),\J\x(t)+\Q_{\CC}-\dot{\bar{\x}}(t)\rangle $ is bounded, and denoted by $ |\phi(t)|\le\Phi(\x_0). $ It is easy to obtain that
\begin{equation}\label{ieq:x-xave}
\begin{aligned}
\|\x(t)-\bar{\x}(t)\|\le&\|\x_0-\bar{\x}(0)\| e^{-K\lambda_{n-1}(\LL_{\G})t}+\sqrt{\frac{\Phi(\x_0)}{2\lambda_{n-1}(\LL_{\G})}}\frac{1}{\sqrt{K}}.
\end{aligned}	
\end{equation}
From  \eqref{def:exp rate to B}, $ \|\x(t)-\bar{\x}(t)\| $ is exponentially convergent to $ B_\delta $ with $$ \delta=\sqrt{\frac{\Phi(\x_0)}{2\lambda_{n-1}(\LL_{\G}))}}\frac{1}{\sqrt{K}} $$ at a rate of $ K\lambda_{n-1}(\LL_{\G}) $. In other words, when $ K $ is large enough, $ \|\x(t)-\bar{\x}(t)\| $ is exponentially convergent to an arbitrarily small neighborhood of the origin at a very fast rate. Moreover, it can be concluded that for any $ \delta>0,  t_0>0, $ there exists $ K_{\delta,t_0} $, such that 
\begin{equation}
\|\x(t)-\bar{\x}(t)\|\le\delta,\quad \forall  K>K_{\delta,t_0}, t>t_0.
\end{equation}
This completes the proof. \hfill $ \Box $

\subsection{Proof\ of\ Theorem\ \ref{thm:convergence}}
Note that $ \dot{\x}=-\J_{\LL}\x+\Q_{\CC} $ has at least one equilibrium point because of Assumption \ref{assumption-exact solution}, and $ \J_{\LL} $ is positive semidefinite because $ K\LL_{\G}\otimes \I_{n^2} $ and $ \J $ are positive semidefinite. Then, for any initial value $ \x_0, $ there exists $ \X^*(\x_0)=\vvec^{-1}(\y^*(\x_0))\in\R^{n\times n} $, such that $ \vvec^{-1}(\x_i(t)) $ converges to $ \X^*(\x_0) $ exponentially, for any $ i\in\V $. Moreover, the rate of the exponential convergence is the minimum non-zero eigenvalue of $ \J_{\LL}, $ denoted by $ r(K). $

\noindent(\romannumeral1) Based on a direct application of the convergence theorem for ``consensus + projection'' flow (Theorem 1, 3 and 5 in\cite{shi2017networkTAC}), we conclude that, for any initial value $ \x_0, $ any $ i\in\V, $
\[
\lim\limits_{t\to \infty}\vvec^{-1}(\x_{i}(t))=\X^*(\x_0)=\vvec^{-1}\Big(\frac{1}{n}\sum_{i=1}^{n}\Pp_{\cap_{i=1}^n\E_i}(\x_i(0))\Big), 
\]
which is a solution to \eqref{eq:Sylvester}.

 \noindent(\romannumeral2) (a) For all $ t\ge0, $
 \begin{equation}\label{eq:thm1_x*}
 \begin{aligned}
 &\sum_{i=1}^n\|\vvec^{-1}(\x_{i}(t))- \X^*(\x_0)\|_F^2=\|\x(t)-\1_n\otimes\y^*(\x_0)\|^2\le (\sum_{i=1}^n\|\x_i(0)-\y^*(\x_0)\|^2)e^{-2r(K)t}.
 \end{aligned}
 \end{equation}
Because $ \J_{\LL} $ is symmetric, the left eigenvector space of its eigenvalue $ r(K) $ is the same as the right eigenvector space, where the base matrix is denoted as $ \Psi_r $. As a special case of 
Lemma \ref{lemma:Eigenvalue Derivatives},  we obtain 
\begin{align*}
\Psi_r^T\frac{\partial \J_{\LL}}{\partial K}\Psi_r=\Psi_r^T(\LL_{\G}\otimes \I_{n^2})\Psi_r,
\end{align*}
which is a positive semidefinite matrix and has eigenvalues in the form of	
$ \partial r(K)/\partial K\ge0 $.
Then, $ r(K) $ is a continuous monotonically non-decreasing function  with respect to $ K. $ 

\noindent(b)	
$ \LL_{\G} $ is a symmetric positive semidefinite matrix with a single zero eigenvalue, which implies that $ \lambda_i(K\LL_{\G}\otimes \I_{n^2}) $ has all non-negative eigenvalues with $ n^2 $ zero.
Denote $$ p:=\rank(\J)=\sum_{i=1}^{n}\rank(\Hh_i^{\dag}\Hh_i)=\sum_{i=1}^{n}\rank(\Hh_i)\le n^2. $$
Following from Lemma \ref{lemma:prop of inverse}, $ \J $ is real symmetric and idempotent, and its eigenvalues can only be zero or one.
Thus, $ 1=\lambda_{1}(\J)=\cdots=\lambda_{p}(\J)>\lambda_{p+1}(\J)=\cdots=\lambda_{n^3}(\J)=0. $
Due to Lemma \ref{lemma:Weyl's inequality},
\begin{align}
\lambda_{n^3-n^2}(\J_{\LL})&\ge
\lambda_{n^3-n^2}(K\LL_{\G}\otimes \I_{n^2})+\lambda_{n^3}(\J)\nonumber=K\lambda_{n-1}(\LL_{\G})>0,\\
\lambda_{n^3-n^2+k}(\J_{\LL})&\le \lambda_{n^3-n^2+k}(K\LL_{\G}\otimes \I_{n^2})+\lambda_{1}(\J)\nonumber=1,\ k=1,\cdots,n^2.\label{eq:r_n3-n2>0,+k<1}
\end{align}
Because $ \mathrm{ker}(\LL_{\G})=\{k\1_n:k\in \R \}, $ \begin{align*}
 &\mathrm{ker}(K\LL_{\G}\otimes \I_{n^2})=\{\col\{\w_1,\cdots,\w_n\}:\w_1=\cdots=\w_n\in \R^{n^2} \},
\end{align*}  and $ \dim(\mathrm{ker}(K\LL_{\G}\otimes \I_{n^2}))=n^2. $
Also, $$\mathrm{ker}(\J_{\LL})=\mathrm{ker}(K\LL_{\G}\otimes \I_{n^2})\cap \mathrm{ker}(\J) $$ because of   positive semidefinite matrices $ K\LL_{\G}\otimes \I_{n^2} $ and $ \J $. Therefore, $ \w:=\col\{\w_1,\cdots,\w_n\}\in \mathrm{ker}(\J_{\LL}) $ if and only if
$ \w_1=\cdots=\w_n $ and \begin{align*}
\J\w&=\sum_{i=1}^{n}\Hh^{\dag}_i\Hh_i\w_i=0;\\
\w_1&=\cdots=\w_n\in \mathrm{ker}(\sum_{i=1}^{n}\Hh^{\dag}_i\Hh_i).
\end{align*}
Hence, from Lemma \ref{lemma:prop of inverse},
\begin{align*}
&\dim(\mathrm{ker}(\J_{\LL}))=\dim(\ker(\sum_{i=1}^{n}\Hh^{\dag}_i\Hh_i)=\dim(\bigcap_{i=1}^n\ker(\Hh^{\dag}_i\Hh_i))= \dim(\bigcap_{i=1}^n\ker(\Hh_i))=\dim(\ker(\Hh)).
\end{align*}  Then 
\begin{align*}
\rank(\J_{\LL})&=n^3-\dim(\ker(\Hh))
=n^3-n^2+\rank(\Hh)>n^3-n^2,\quad \text{if}\ \rank(\Hh)\neq0.
\end{align*}
Now we can conclude that, $r(K)=\lambda_{\rank(\J_{\LL})}(\J_L)=\lambda_{n^3-n^2+\rank(\Hh)}(\J_{\LL})\le1 $, is always bounded.

\noindent(c)
It has been proved that $ r(K) $ is always upper bounded and $ r(K)\le1 $. Since also $ r(K) $ increases with increasing $ K, $ there must exist a limit $ r^*=\lim\limits_{K\to \infty}r(K) $.
To prove $ r^*=\lambda_{\rank(\Hh)}((\sum\nolimits_{i=1}^{n}\Hh^{\dag}_i\Hh_i)/n), $ we take four steps. For convenience below, we define $ r_0=\lambda_{\rank(\Hh)}((\sum\nolimits_{i=1}^{n}\Hh^{\dag}_i\Hh_i)/n). $

\noindent{\bf Step 1:} In this step, we prove $$ \|\x_{ave}(t)-\y^*(\x_0)\|\le c_1e^{-r(K)t} $$ for some $ c_1>0. $ 
Combining $ \x_{ave}(t)=\frac{1}{n}\sum_{i=1}^n\x_i(t) $ and \eqref{eq:thm1_x*}, we get the convergence of $ \x_{ave}(t), $
\begin{equation}\label{eq:xave_y*}
\begin{aligned}
&\|\x_{ave}(t)-\y^*(\x_0)\|^2
=\frac{1}{n^2}\|\sum_{i=1}^n\x_{i}(t)-n\y^*(\x_0)\|^2
\le\frac{1}{n^2}\sum_{i=1}^n\|\x_i(t)-\y^*(\x_0)\|^2\le c_1^2(\x_0)e^{-2r(K)t},
\end{aligned}
\end{equation}
where $ c_1(\x_0)=1/n\|\x_0-\1_n\otimes\y^*(\x_0)\| $ and $ c_1(\x_0) $ is denoted as $ c_1 $ for  simplicity.
Then $ \x_{ave}(t) $ converges to $ \y^*(\x_0) $ exponentially at the rate $ r(K). $ 

\noindent{\bf Step 2:} In this step, we prove $ r^*\le r_0 $. Summing equations in \eqref{eq:Flow-n-nodes} from $ i=1 $ to $ i=n $, we obtain 
$ \sum_{i=1}^{n}\dot{\x}_i(t)=\sum_{i=1}^{n}(\Pp_{\E_i}(\x_i(t))-\x_i(t)). $
Next, there holds 
\begin{equation}\label{eq:x_ave}
\begin{aligned}
\dot{\x}_{ave}(t)=&\frac{1}{n}\sum_{i=1}^{n}(\Pp_{\E_i}(\x_i(t))-\x_i(t))\\=&\frac{1}{n}\sum_{i=1}^{n}(-\Hh^{\dag}_i\Hh_i\x_{i}(t)+\Hh^{\dag}_i\CC_i)\\
=&-\frac{1}{n}(\sum_{i=1}^{n}\Hh^{\dag}_i\Hh_i)\x_{ave}(t)+\frac{1}{n}\sum_{i=1}^{n}\Hh^{\dag}_i\CC_i+\frac{1}{n}\sum_{i=1}^{n}(\Hh^{\dag}_i\Hh_i(\x_{ave}(t)-x_i(t))).
\end{aligned}
\end{equation}
Denoting $ \sigma=\max_{i\in\{1,\cdots,n\}}\lambda_1(\Hh^{\dag}_i\Hh_i), $ we rewrite \eqref{eq:x_ave} as
\begin{align*}
\frac{d}{dt}(\x_{ave}(t)-\y^*(\x_0))=&-\frac{1}{n}(\sum_{i=1}^{n}\Hh^{\dag}_i\Hh_i)(\x_{ave}(t)-\y^*(\x_0))+\underbrace{\frac{1}{n}\sum_{i=1}^{n}(\Hh^{\dag}_i\Hh_i(\x_{ave}(t)-x_i(t)))}_{\nu(t)}.
\end{align*}
 We have $ \|\nu(t)\|=O(e^{-r(K)t}) $ because  $$ \|\Hh^{\dag}_i\Hh_i(\x_{ave}(t)-\x_i(t))\|=O(e^{-r(K)t}), $$ which is owing to
\begin{itemize}
	\item $	
	 \|\Hh^{\dag}_i\Hh_i(\x_{ave}(t)-\x_i(t))\|\le\sigma\|\x_{ave}(t)-\x_i(t)\| \le\sigma\|\x_{ave}(t)-\y^*(\x_0
	)\|+\sigma\|\y^*(\x_0)-\x_i(t)\|;$
	\item $ \|\x_{ave}(t)-\y^*(\x_0)\|=O(e^{-r(K)t})$ from \eqref{eq:xave_y*}; 
	\item $ \|\y^*(\x_0)-\x_i(t)\|=O(e^{-r(K)t}) $ from \eqref{eq:thm1_x*}.
\end{itemize}

Now we prove that $ r^*\le r_0 $ by contradiction. Suppose $ r^*>r_0, $ there exist $\epsilon>0$ and $ K_{\epsilon}>0 $, such that $ r(K_{\epsilon})=r_0+\epsilon< r^* $ and $ \|\x_{ave}(t)-\y^*(\x_0)\|=O(e^{-(r_0+\epsilon)t}) $ due to \eqref{eq:xave_y*}. Considering Lemma \ref{lemma:LTIorder} and $ \|\nu(t)\|=O(e^{-r(K_\epsilon)t})=O(e^{-(r_0+\epsilon)t}) $, we get $ \|\x_{ave}(t)-\y^*(\x_0)\|=\Theta(e^{-r_0t}) $, which leads to a contradiction.

\noindent	{\bf Step 3:} In this step, we prove $$ \|\x_{ave}(t)-\y^*(\x_0)\|\le c_2e^{-r_0t}+\mu $$ for some $ c_2>0 $ and any $ \mu>0 $.
From the elementary inequality $$ \frac{\sum_{i=1}^nx_i}{n}\le\sqrt{\frac{\sum_{i=1}^{n}x_i^2}{n}} $$ and Lemma \ref{lemma:step3part}, for any $ \delta, t_0>0, $ there exists $ K_{\delta,t_0} $, such that \begin{align*}\label{ieq:sumx_ave-xi}
&\sum_{i=1}^n\|\x_{ave}(t)-\x_i(t)\|\le (n\sum_{i=1}^n\|\x_{ave}(t)-\x_i(t)\|^2)^{\frac{1}{2}}=\left(n\|\x(t)-\bar{\x}(t)\|^2\right)^{\frac{1}{2}}\le\sqrt{n}\delta,\quad \forall  K>K_{\delta,t_0}, t>t_0.
\end{align*}		 	
Due to $ \|\x_{ave}(t)-\y^*(\x_0)\|<c_1e^{-r(K)t} $ in {\bf Step 1} and  Lemma \ref{lemma:step3part},  for any $ \delta, t_0>0, $ any $ K>K_{\delta,t_0}, t>t_0 $,  denoting   $ \omega(t)=2\big\langle\x_{ave}(t)-\y^*(\x_0),\nu(t)\big\rangle $, we have	
\begin{equation}
\begin{aligned}
|\omega(t)|=&|\frac{2}{n}(\x_{ave}(t)-\y^*(\x_0))^T(\sum_{i=1}^{n}(\Hh^{\dag}_i\Hh_i(\x_{ave}(t)-\x_i(t)))|\\\le&\frac{2}{n}\|\x_{ave}(t)-\y^*(\x_0)\| \sigma(\sum_{i=1}^{n}\|\x_{ave}(t)-\x_i(t)\|)
\le\frac{2\sigma}{n} c_1e^{-r(K)t}\sqrt{n}\delta\le 
(\frac{2\sigma}{\sqrt{n}} c_1)\delta.
\end{aligned}
\end{equation}
Thus, we obtain $ |\omega(t)|\to 0 $ as $\delta\to 0, $ for any $ K>K_{\delta,t_0}, t>t_0 $. Because of the arbitrariness of $ t_0, $ it easy to prove that $ \omega(t) $ is always bounded for all $ t>0. $
Then we think about $ \|\x_{ave}(t)-\y^*(\x_0)\|^2, $
\begin{equation}\label{ieq:dx-y2}
\begin{aligned}
&\frac{d}{dt}\|\x_{ave}(t)-\y^*(\x_0)\|^2
\le-2r_0\|\x_{ave}(t)-\y^*(\x_0)\|^2+\omega(t).
\end{aligned}
\end{equation}
With the help of the Gr$ \ddot{\mathrm{o}} $nwall Inequality, we obtain 
\begin{equation}\label{Ieq:|x_ave-y*|^2}
\begin{aligned}
\|\x_{ave}(t)-\y^*(\x_0)\|^2\le& \|\x_{ave}(0)-\y^*(\x_0)\|^2e^{-2r_0t}+\int_0^te^{-2r_0(t-s)}\omega(s)ds.
\end{aligned}\end{equation}
Thus,  with $ \beta_{t_0}:=\max_{0<t
	\le t_0}|\omega(t)| $,  for  $ t>0, $
\begin{equation*}
\begin{aligned}
\int_0^te^{-2r_0(t-s)}\omega(s)ds&<\frac{\beta_{t_0} e^{2r_0t_0}}{2r_0}e^{-2r_0t}+\frac{\sigma c_1\delta}{r_0\sqrt{n}},~t>t_0;\\
\int_0^te^{-2r_0(t-s)}\omega(s)ds&\le\int_0^{t_0}e^{-2r_0(t-s)}\omega(s)ds
<\frac{\beta_{t_0} e^{2r_0t_0}}{2r_0}e^{-2r_0t},~t\le t_0.
\end{aligned}
\end{equation*}
That is, for all $ t>0,$ \[\int_{0}^te^{-2r_0(t-s)}\omega(s)ds<\frac{\beta_{t_0} e^{2r_0t_0}}{2r_0}e^{-2r_0t}+\frac{\sigma c_1\delta}{r_0\sqrt{n}}. \] Specifically, setting $ t_0=1, $ the inequality \eqref{Ieq:|x_ave-y*|^2} implies that, for all $ t>0 $,
\begin{align*}\label{Ieq:|x_ave-y*|^2again}
\|\x_{ave}(t)-\y^*(\x_0)\|^2\le& (\|\x_{ave}(0)-\y^*(\x_0)\|^2+\frac{\beta_1 e^{2r_0}}{2r_0})e^{-2r_0t}+\frac{\sigma c_1\delta}{r_0\sqrt{n}}.
\end{align*}
Then, for any $ \mu>0 $ with $ \delta=r_0\sqrt{n}\mu^2/(\sigma c_1), $ for all $ K>K_{\delta,1}, $ there holds 
\begin{align*}
 \|\x_{ave}(t)-\y^*(\x_0)\|^2\le c_2^2(\x_0)e^{-2r_0t}+\mu^2,\quad
 \text{with}~ c_2:=c_2^2(\x_0)=\|\x_{ave}(0)-\y^*(\x_0)\|^2+\frac{\beta_1 e^{2r_0}}{2r_0}.
\end{align*} 
Further,  for all $ K>K_{\delta,1}, t>0 $,
\begin{equation}
\|\x_{ave}(t)-\y^*(\x_0)\|\le(c_2^2e^{-2r_0t}+\mu^2)^{\frac{1}{2}}<c_2e^{-r_0t}+\mu.
\end{equation} 

\noindent  {\bf Step 4:} Let us complete the proof of $ r_0=r^* $, while
 $ r_0\ge r^* $ has been shown  in {\bf Step 2}. We now prove $ r_0\le r^* $ by contradiction.
Assume $ r_0>r^*=\sup_{K>0}r(K). $ Then there exists $ \eta>0 $ satisfying $ r_0>r^*+\eta, $ such that for all $ K\ge K_{\delta,1}, t>0, $
$$\|\x_{ave}(t)-\y^*(\x_0)\|<c_2e^{-r_0t}+\mu<c_2e^{-(r^*+\eta)t}+\mu.$$
Following from \eqref{ieq:x-xave}, there holds,  
\begin{equation}
\begin{aligned}
& \|\x(t)-\1_n\otimes\y^*(\x_0)\|\\=& \|\x(t)-\bar{\x}(t)+\1_n\otimes(\x_{ave}(t)-\y^*(\x_0))\|\\
\le& \|\x_0-\bar{\x}(0)\| e^{-K\lambda_{n-1}(\LL_{\G})t}+\sqrt{\frac{\Phi(\x_0)}{2\lambda_{n-1}(\LL_{\G})}}\frac{1}{\sqrt{K}}+n(c_2e^{-(r^*+\eta)t}+\mu).
\end{aligned}
\end{equation}
When $ K_m=\max\{K_{\delta,1},\frac{r^*+\eta}{\lambda_{n-1}(\LL_{\G})}\} $, for all $ K>K_m, $ we have
\begin{equation}
\|\x(t)-\1_n\otimes\y^*(\x_0)\|
\le c_3e^{-(r^*+\eta)t}+g(K,\mu),\quad t>0,
\end{equation}
where $ c_3=\|\x_0-\bar{\x}(0)\|+nc_2 $ and $$  g(K,\mu)=\sqrt{\frac{\Phi(\x_0)}{2\lambda_{n-1}(\LL_{\G})}}\frac{1}{\sqrt{K}}+n\mu. $$	
By Lemma \ref{lemma:LTIorder}, $ \|\x(t)-\1_n\otimes\y^*(\x_0)\|=\Theta(e^{-r(K)t}). $ There is $ t'>0 $ and a positive constant $ p(\x_0) $ depending on $ \x_0 $, such that, for all $ K>0, t>t',$ $$ \|\x(t)-\1_n\otimes\y^*(\x_0)\|\ge p(\x_0)e^{-r(K)t}. $$
Then, for any $ \mu>0 $, for all $ K\ge K_m, t>t' $, there holds
\begin{equation}
\begin{aligned}
p(\x_0)e^{-r^*t}&\le p(\x_0)e^{-r(K)t}\le\|\x(t)-\1_n\otimes\y^*(\x_0)\|<c_3e^{-(r^*+\eta)t}+g(K,\mu).
\end{aligned}
\end{equation}
Equivalently, $ p(\x_0)e^{-r^*t}<c_2e^{-(r^*+\eta)t}+g(K,\mu) $ for any $ \mu>0 $ and all $ K\ge K_m, t>t'. $ However, the positive term $ g(K,\mu) $ can be arbitrarily small with $ K $ large enough and $ \mu $ small enough, which leads to a contradiction.  
Therefore, $ r^*=r_0=  \lambda_{\rank(\Hh)}(\frac{1}{n}(\sum_{i=1}^{n}\Hh^{\dag}_i\Hh_i)). $ The proof  has been completed. 
$ \hfill $ \qedsymbol

\section{Proof of Theorem \ref{thm:sym-convergence}}\label{App:sym-con}
 Denote by $ \Ss^n $ the set of all $ n\times n $ real symmetric matrices. Then $ \Ss^n $ is convex.
The projector onto $ \Ss^n $, $ \Pp_{\Ss^n}(\cdot):\R^{n\times n}\to\Ss^n $ is an orthogonal projection with concrete expression $ \Pp_{\Ss^n}(\X)=(\X+\X^T)/2. $ In order to vectorize it, we define $ \Ss_{nn}=\{\y\in\R^{n^2}:\y=\vvec(\X), \ \mathrm{for\ some\ } \X\in\Ss^n\}. $ As a result, a projection mapping $ \Pp_{\Ss_{nn}} $ satisfies
\begin{align*} &\Pp_{\Ss_{nn}}(\y)=\vvec\Big(\Pp_{\Ss^n}(\vvec^{-1}(\y))\Big)
=\vvec\left(\frac{1}{2}(\vvec^{-1}(\y)+(\vvec^{-1}(\y))^T\right)
=\frac{1}{2}(\y+\mathbf{P}_{n^2}\y),
\end{align*}
where $ \mathbf{P}_{n^2}\in\R^{n^2} $ is an elementary matrix obtained by swapping row $ (k-1)n+j $ and row $ (j-1)n+k $  for every $ k=1,\cdots,n $ and $ k<j\le n $ of the identity matrix.
Then the flow \eqref{eq:Flow-n-nodes and symmetric} can be represented as 
\begin{equation}\label{eq:Flow-n-nodes and symmetric'}
\begin{aligned}
\dot{\x}_i
=&K\Big(\sum_{j\in \N_i}(\x_j-\x_i)\Big)-\Hh_i^{\dag}\Hh_i\x_i+\Hh_i^{\dag}\CC_i+\frac{K_s}{2}(\mathbf{P}_{n^2}-\I_{n^2})\x_i,
\quad i\in\V.
\end{aligned}
\end{equation} 
The compact form of \eqref{eq:Flow-n-nodes and symmetric'} is
\begin{equation}\label{eq:Flow-X and Symmetric}
\begin{aligned}
\dot{\x}=-\underbrace{(K\LL_{\G}\otimes \I_{n^2}+\J_p)}_{\J_{LP}}\x+\Q_{\CC},
\end{aligned}
\end{equation}
where $ \LL_{\G} $ is the Laplacian matrix, 
\begin{align*}
\J_p&=\mathrm{diag}\{\Hh^{\dag}_i\Hh_i+\frac{K_s}{2}(\I_{n^2}-\mathbf{P}_{n^2}), i\in\V\}\in \mathbb{R}^{n^3\times n^3},\\
\Q_{\CC}&=\col\{\Hh^{\dag}_1\CC_1, \cdots, \Hh^{\dag}_n\CC_n\}.
\end{align*} Besides, we know that both $ \Hh^{\dag}_i\Hh_i $ and $ \frac{1}{2}\I_n\otimes(\I_{n^2}-\mathbf{P}_{n^2}) $ only have eigenvalues $ 1 $ and $ 0. $

Note that system \eqref{eq:Flow-X and Symmetric} has at least one equilibrium point and $ (\cap_{i=1}^n\E_i)\bigcap\Ss_{nn}\ne\emptyset $  because of the existence of a symmetric solution. Since the matrix $ \J_{LP} $ is positive semidefinite, we conclude that, for any initial value $ \x_0 $, there exists $ \X_s^*(\x_0)\in\Ss^{n\times n} $, such that $ \vvec^{-1}(\x_i(t)) $ along the flow \eqref{eq:Flow-n-nodes and symmetric'}  converges to $ \X_s^*(\x_0) $ exponentially, for any $ i\in\V $. Moreover, the rate of the exponential convergence is the minimum non-zero eigenvalue of $ \J_{LP}, $ denoted by $ r_s(K,K_s)=\min\{\lambda:\lambda\in\mathrm{spec}(\J_{LP}), \lambda\ne0\}. $

Because $ \rank(\Hh_i^{\dag}\Hh_i)=\rank(\Hh_i) $ and $ \rank(\I_{n^2}-\mathbf{P}_{n^2})=n(n-1)/2, $ we have
\begin{align*}
\rank(\J_p)\le& n\big(\rank(\Hh_i^{\dag}\Hh_i)+\rank(\frac{K_s}{2}(\I_{n^2}-\mathbf{P}_{n^2}))\big)\le n(n+\frac{n^2-n}{2})=\frac{n^3+n^2}{2}<n^3,\ \text{if}\ n\ge2.
\end{align*} 
What's more, $ \lambda_{n^3}(\J_p)=0 $, $ \lambda_1(\J_p)\le1+K_s $ following from \begin{align*}
\lambda_1(\J_p)\le&\lambda_1(\diag\{\Hh_i^{\dag}\Hh_i, i\in\V\})+\lambda_1(\frac{K_s}{2}\I_n\otimes(\I_{n^2}-\mathbf{P}_{n^2})).
\end{align*} Therefore, if $ \rank(\Hh)\ne0 $,
\begin{align*}
\rank(\J_{LP})&=n^3-\dim(\ker(\J_{LP}))\ge n^3-(n^2-\text{rank}(\Hh))>n^3-n^2,  
\end{align*}	since
\begin{align*}
\dim(\ker(\J_{LP}))=&\dim\Big((\cap_{i=1}^n\ker(\Hh_i^{\dag}\Hh_i))\cap\ker(\frac{K_s}{2}(\I_{n^2}-\mathbf{P}_{n^2}))\Big)\\=& \dim\Big((\cap_{i=1}^n\ker(\Hh_i))\cap\ker(\I_{n^2}-\mathbf{P}_{n^2})\Big)\\
\le&\min\{\dim(\ker(\Hh)),\dim(\ker(\I_{n^2}-\mathbf{P}_{n^2}))\}\\
\le&\min\{n^2-\text{rank}(\Hh),\frac{n^2+n}{2}\}.
\end{align*}
Now we can conclude that, for all $ K>0, $ 
\begin{align*}
&\lambda_{n^3-n^2+k}(\J_{LP})\le\lambda_{n^3-n^2+k}(K\LL_{\G}\otimes \I_{n^2})+\lambda_1(\J_p)=0+\lambda_1(\J_p)\le1+K_s,  \text{\ for\ all} \ k=1,\cdots,n^2. 
\end{align*}
 Similarly, we also have
\begin{align*}
\lambda_{n^3-n^2+k}(\J_{LP})\le&\lambda_{n^3-n^2+k}(\mathrm{diag}\{\frac{K_s}{2}(\I_{n^2}-\mathbf{P}_{n^2}), i\in\V\})+\lambda_{1}(K\LL_{\G}\otimes \I_{n^2})+\lambda_1(\diag\{\Hh_i^{\dag}\Hh_i, i\in\V\}) \\
=&0+K\lambda_{1}(\LL_{\G})+1\quad \text{for all} \ k=1,\cdots,n^2, 
\end{align*} due to $$ \rank(\frac{K_s}{2}(\I_{n^2}-\mathbf{P}_{n^2}), i\in\V\})=(n^3-n^2)/2\le n^3-n^2. $$ As a result, $$ r_s(K,K_s)\le\min\{1+K_s, 1+K\lambda_1(\LL_{\G})\} $$ for all $ K, K_s>0. $ The proof of Theorem \ref{thm:sym-convergence} has been completed.
$ \hfill $\qedsymbol

\section{Proof of Theorem \ref{thm:row+column convergence}}\label{App:full}
 Under the $ \A $ Row/$ \B $-$ \CC $ Column Partition, the equation   \eqref{eq:Sylvester} is equivalent to (\cite{zeng2018distributedlmeTAC})
 \begin{equation}\label{eq:equivalent_w}
 \begin{aligned}
 &\ee_i(\A^T)_i^T\X+\X\B_i\ee_i^T-((\LL_{\G}^T)_i^T\otimes\I_{n})\mathbf{Z}=\CC_{i}\ee_i^T, \ i\in\V, \ \text{for\ some}\ \mathbf{Z}\in\R^{n^2\times n},
 \end{aligned} 
 \end{equation}
 where $ (\LL_{\G}^T)_i^T $ represents the $ i $-th row of the Laplacian matrix $ \LL_{\G} $.
 Taking advantage of Kronecker product, we rewrite \eqref{eq:equivalent_w}, for all $ i\in\V, $
 \begin{equation}\label{eq:vec_equivalent_w}
 \begin{aligned}
 \vvec(\CC_{i}\ee_i^T)=&\big(\I_n\otimes(\ee_i(\A^T)_{i}^T)+(\ee_i\B_{i}^T)\otimes\I_n\big)\vvec(\X)-\I_n\otimes((\LL_{\G}^T)_i^T\otimes\I_n)\vvec(\mathbf{Z})\\=& [\I_n\otimes(\ee_i(\A^T)_{i}^T)+(\ee_i\B_{i}^T)\otimes\I_n, -\I_n\otimes((\LL_{\G}^T)_i^T\otimes\I_n)]\y
 \end{aligned}\end{equation}with  $ \y=\vvec([\X,\mathbf{Z}]). $
  Therefore, 
 the flow \eqref{eq:row+Flow-n-nodes} is an extended form of the ``consensus $ + $ projection'' flow \eqref{eq:Flow-n-nodes} with respect to the augmented variable $ \y_i=\col\{\x_i,\z_i\} $.
 
 	Since Assumption \ref{assumption-exact solution} holds,  \eqref{eq:equivalent_w} has at least one solution. We replace $ \x_{i}(t) $ in Theorem \ref{thm:convergence} with $ \y_i(t) $ and then conclude that $ \y_i(t) $ along the flow \eqref{eq:row+Flow-n-nodes} converges to $ 1/n\sum_{i=1}^n\Pp_{\cap_{i=1}^d\E_i^{\text{Aug}}}(\y_i(0)) $ exponentially, which is a solution to \eqref{eq:vec_equivalent_w}. Moreover,  we can conclude that, there exists $ \X^*(\y_0)\in\R^{n\times n}, $ such that  
 	$ \vvec^{-1}(\x_i(t)) $  converges to $ \X^*(\y_0) $ exponentially, for all $ i\in\V $. 
 	$ \hfill $\qedsymbol

\section{Proof of  Theorem \ref{thm:doublecase}}\label{App:double-con}
\subsection{Key Lemma}
Rewrite \eqref{eq:double-flow 1st partition_x} as a compact form:
\begin{align}\label{eq:x_z 1st_matrix}
\begin{bmatrix}
\dot{\x}\\
\dot{\z}
\end{bmatrix}=-\Gg\begin{bmatrix}
\x\\
\z
\end{bmatrix}+\begin{bmatrix}
\bar{\M}^T\bar{\CC}\\
-\bar{\CC}
\end{bmatrix}.
\end{align}
Denote $ \M_{\LL}:=[
\bar{\M}, -\bar{\LL}], $ and then
\begin{equation}\label{eq:Gbar}
\begin{aligned}
\bar{\Gg}:=&\begin{bmatrix}
\I_{n^3} & \0\\
\0 & \bar{\LL}
\end{bmatrix}\Gg=\M_{\LL}^T\M_{\LL}+\diag\{K(\LL_{\G}\otimes\I_{n^2}),\0_{n^3\times n^3}\}.
\end{aligned}
\end{equation}
So $ \bar{\Gg} $ is positive semi-definite as the sum of two positive semi-definite matrices.
The following lemma is a generalization of Lemma \ref{lemma:G_eigenvalue>0or=0}.

\begin{lem}\label{lemma:G-Gbar-prop}
	Suppose that
	$ \Gg $ and $ \bar{\Gg} $ are defined as in \eqref{eq:x_z 1st_matrix} and \eqref{eq:Gbar}. Then they obtain following properties.
	\begin{enumerate}
		\item[(\romannumeral1)] $ \Gg $ is a non-defective matrix with all eigenvalues being real non-negative;
		\item[(\romannumeral2)] $ \ker(\Gg)=\ker(\bar{\Gg}) $ and $ \rank(\Gg)=\rank(\bar{\Gg}). $
	\end{enumerate}
\end{lem}
\noindent$ Proof. $
(\romannumeral1) Using  Lemma \ref{lemma:G_eigenvalue>0or=0} and the fact that the matrices $ K(\LL_{\G}\otimes\I_{n^2}) $ and $ \bar{\LL} $ are positive semi-definite, we conclude that all eigenvalues of $ \Gg $ are greater than or equal to 0. Moreover, the possible zero eigenvalue must be non-defective. We prove that any positive eigenvalue $ \lambda $ is non-defective by contradiction. Suppose that $ \lambda>0 $ is defective, namely, there exists a non-zero vector $ \col\{\vv_1,\vv_2\} $ such that 
\begin{equation}\label{q5}
(\Gg-\lambda \I_{2n^3})\begin{bmatrix}
\vv_1\\\vv_2\end{bmatrix}=\begin{bmatrix}
\uu_1\\\uu_2\end{bmatrix}\neq 0,\ (\Gg-\lambda \I_{2n^3})\begin{bmatrix}
\uu_1\\\uu_2\end{bmatrix}=0.
\end{equation}	
Premultiplying the first equation in \eqref{q5} by $ [\uu_1^T,\uu_2^T]\begin{bmatrix}
\I_{n^3} & \0\\
\0 & \bar{\LL}
\end{bmatrix} $, we have
\begin{equation}\label{q6'}
[\uu_1^T,\uu_2^T](\bar{\Gg}-\lambda\begin{bmatrix}
\I_{n^3} & \0\\
\0 & \bar{\LL}
\end{bmatrix})\begin{bmatrix}\vv_1\\\vv_2\end{bmatrix}=\uu_1^T\uu_1+\uu_2^T\bar{\LL}\uu_2.
\end{equation}
For the second equation in \eqref{q5},
\begin{equation}\label{q7}
\begin{bmatrix}
\I_{n^3} & \0\\
\0 & \bar{\LL}
\end{bmatrix}(\Gg-\lambda \I_{2n^3})\begin{bmatrix}
\uu_1\\\uu_2\end{bmatrix}=(\bar{\Gg}-\lambda\begin{bmatrix}
\I_{n^3} & \0\\
\0 & \bar{\LL}
\end{bmatrix})\begin{bmatrix}\uu_1\\\uu_2\end{bmatrix}=0.
\end{equation}		
Because $ \bar{\Gg} $ and $ \bar{\LL} $ are symmetric, substituting the transpose of  \eqref{q7} into \eqref{q6'} yields 
$ \uu_1^T\uu_1+\uu_2^T\bar{\LL}\uu_2=0. $
Because of    positive semi-definite matrix $ \bar{\LL} $, there hold
\begin{equation}\label{q9}
\uu_1=\0, \quad \bar{\LL}\uu_2=\0.
\end{equation}
Rewriting the second equation in \eqref{q5} according to the definition of $ \Gg $, we have
\begin{equation}\label{q10}
\begin{cases}
\left(\bar{\M}^T\bar{\M}+K(\LL_{\G}\otimes\I_{n^2})-\lambda\I_{n^3}\right)\uu_1-\bar{\M}^T\bar{\LL}\uu_2=\0,\\
-\bar{\M}\uu_1+(\bar{\LL}-\lambda\I_{n^3})\uu_2=\0.
\end{cases}
\end{equation}
Substitution of \eqref{q9} into \eqref{q10} yields $ \lambda\I_{n^3}\uu_2=\0, $ then $ \uu_2=\0 $ since $ \lambda>0. $ Clearly, $ \col\{\uu_1,\uu_2\}=\0 $ leads to a contradiction. Thus, any non-zero eigenvalue $ \lambda $ of $ \Gg $ is non-defective and $ \Gg $ is a non-defective matrix.

\noindent(\romannumeral2) Note that $ \ker(\Gg)\subset\ker(\bar{\Gg}) $ follows from the definition of $ \bar{\Gg} $. On the other hand, with any $ 
\col\{\bar{\vv}_1,\bar{\vv}_2\}\in \ker(\bar{\Gg}) $,  \eqref{eq:Gbar} yields
\begin{equation}\label{q11}
\begin{cases}
\bar{\M}\bar{\vv}_1-\bar{\LL}\bar{\vv}_2=\0,\\
K(\LL_{\G}\otimes\I_{n^2})\bar{\vv}_1=\0.
\end{cases}
\end{equation}
  $ \Gg\col\{\bar{\vv}_1,\bar{\vv}_2\}=\0 $ holds as a consequence of  \eqref{q11}. Thus, $ \col\{\bar{\vv}_1,\bar{\vv}_2\}\in \ker(\Gg) $ and  $ \ker(\bar{\Gg})\subset\ker(\Gg). $ As a result, $ \ker(\bar{\Gg})=\ker(\Gg), $ and moreover $ \rank(\Gg)=\rank(\bar{\Gg}). $ \hfill$ \Box $
\subsection{Proof of  Theorem \ref{thm:doublecase}}
 The convergence of the flow \eqref{eq:double-flow 1st partition_xi_j} as a direct application of Theorem 1 in \cite{wang2019scalable}. Now we prove the properties of the exponential convergence rate $ r^*(K) $.
Following from Lemma \ref{lemma:G-Gbar-prop}, $ \Gg $ is a non-defective matrix with non-negative eigenvalues, then the rate of exponential convergence is $ r^*(K)=\min\{\lambda\in\mathrm{spec}(\Gg),\lambda\ne0 \}, $ namely, $ r^*(K)=\lambda_{\rank(\Gg)}(\Gg). $ 
Assume that $ r^*(K) $ is a eigenvalue of $ \Gg $ with multiplicity $ k_r $ and $ \col\{\mathbf{\varphi}_1, \mathbf{\varphi}_2\} $ (with $ \mathbf{\varphi}_1, \mathbf{\varphi}_2\in \R^{n^3} $) is a right eigenvector associated with $ r^*(K) $. Then $ \col\{\mathbf{\varphi}_1, \bar{\LL}\mathbf{\varphi}_2\} $ is a corresponding left eigenvector, which can be proved by direct calculation and the fact $ \bar{\LL}(\bar{\LL}-r^*(K)\I_{n^3})=(\bar{\LL}-r^*(K)\I_{n^3})\bar{\LL} $. Specifically, denote $ \col\{\mathbf{\theta}_1, \mathbf{\theta}_2\}:=\col\{\mathbf{\varphi}_1, \bar{\LL}\mathbf{\varphi}_2\} $ and
\begin{align*}
\Gg\begin{bmatrix}
\mathbf{\varphi}_1\\\mathbf{\varphi}_2
\end{bmatrix}=r^*(K)\begin{bmatrix}
\mathbf{\varphi}_1\\\mathbf{\varphi}_2
\end{bmatrix}, \qquad
\begin{bmatrix}
\mathbf{\theta}_1\\\mathbf{\theta}_2
\end{bmatrix}^T\Gg=r^*(K)\begin{bmatrix}
\mathbf{\theta}_1\\\mathbf{\theta}_2
\end{bmatrix}^T.
\end{align*}
Thus, we  find two base matrices $ \Psi_r=\col\{\Psi_{r1}, \Psi_{r2}\} $ and $ \Theta_r=\col\{\Theta_{r1},\Theta_{r2}\}=\col\{\Psi_{r1}, \bar{\LL}\Psi_{r2}\} $  as the right and the left eigenvector space, respectively, and there holds $ \Theta_r^T\Psi_r=\I_{k_r} $. 
Consider the matrix 
\begin{align*}\label{derivative matrix}
\Theta^T_r\frac{\partial \Gg}{\partial K}\Psi_r&=[\Theta^T_{r1},\Theta^T_{r2}]\begin{bmatrix}
\LL_{\G}\otimes \I_{n^2} & \0_{n^3\times n^3}\\
\0_{n^3\times n^3} & \0_{n^3\times n^3}
\end{bmatrix}
\begin{bmatrix}
\Psi_{r1}\\\Psi_{r2}
\end{bmatrix}\\&=\Theta^T_{r1}(\LL_{\G}\otimes \I_{n^2})\Psi_{r1}
=\Psi_{r1}^T(\LL_{\G}\otimes \I_{n^2})\Psi_{r1},
\end{align*}
which is a positive semidefinite matrix. By Lemma \ref{lemma:Eigenvalue Derivatives},
$ \partial r^*(K)/\partial K\ge0 $. 
As a result, $ r^*(K) $ is a monotonically non-decreasing function with respect to $ K $. 

Next
from \eqref{eq:Gbar}, we obtain \[ \ker(\bar{\Gg})\supseteq\begin{bmatrix}
\ker\bar{\M}\cap\ker(\LL_{\G}\otimes\I_{n^2})\\\ker(\bar{\LL})
\end{bmatrix}=\begin{bmatrix}
\cap_{i=1}^n\ker(\M_i)\\\ker(\bar{\LL})
\end{bmatrix}. \]
Then due to Lemma \ref{lemma:G-Gbar-prop}, 
\begin{align*}
k=&\rank(\Gg)=\rank(\bar{\Gg})=2n^3-\dim(\ker(\bar{\Gg}))\le2n^3-n^2-\dim(\cap_{i=1}^n\ker(\M_i)),
\end{align*} 
which implies the conclusion. \hfill$ \Box $

\end{document}